\documentclass[11pt]{amsart}
\usepackage{geometry}   
\usepackage{pinlabel}             
\usepackage{graphicx}
\usepackage{amssymb, amsthm}
\usepackage{epstopdf}
 \newtheorem{theorem}{Theorem}[section]
    \newtheorem{corollary}[theorem]{Corollary}
   \newtheorem{lemma}[theorem]{Lemma}
    \newtheorem{proposition}[theorem]{Proposition}  
       
    \theoremstyle{definition}
\newtheorem{definition}[theorem]{Definition}
\newtheorem{remark}[theorem]{Remark}
\newtheorem{example}[theorem]{Example}
\newtheorem{examples}[theorem]{Examples}

\newtheorem*{notation}{Notation}

  \newcommand{\Z}{\ensuremath{{\mathbb{Z}}}}
\newcommand{\R}{\ensuremath{{\mathbb{R}}}}
  \newcommand{\G}{\Gamma}
  \newcommand{\iso}{\cong}

\newcommand{\Fone}{F\langle U \rangle}         
\newcommand{\Ftwo}{F\langle W \rangle}        
\newcommand{\AG}{A_\G}            
\newcommand{\OG}{\ensuremath{{\mathcal{O}(A_\Gamma)}}} 
\newcommand{\X}{\ensuremath{{\mathfrak{X}}}}                        
\newcommand{\p}{^\perp}
\newcommand{\TvP}{T_v^\perp}           
\newcommand{\TwP}{T_w^\perp}
\newcommand{\TuP}{T_u^\perp}

\newcommand{\lv}{\ell(v)}          
\newcommand{\tv}{\tau(v)}          
\newcommand{\bv}{\beta(v)}
\newcommand{\dv}{\delta(v)}
\newcommand{\dov}{\delta_0(v)}
\newcommand{\cdv}{\delta_C(v)}
  
\newcommand{\kerR}{K_0}  
\newcommand{\kerP}{K_P}   

\title{Automorphisms of 2-dimensional right-angled Artin groups}
\author{Ruth Charney, John Crisp, and Karen Vogtmann}
\date{}  

\begin{document}

\begin{abstract}      We study the outer automorphism group of a right-angled Artin group $\AG$ in the case where the defining graph $\G$ is connected and triangle-free.   We give an algebraic description of $Out(\AG)$ in terms of  {\it maximal join\,} subgraphs in $\G$ and prove that the Tits' alternative holds for $Out(\AG)$.  We construct an analogue of {\it outer space\,} for $Out(\AG)$ and prove that it is finite dimensional, contractible, and has a proper action of $Out(\AG)$.  We show that $Out(\AG)$ has finite virtual cohomological dimension, give upper and lower bounds on this dimension and construct a spine for outer space realizing the most general upper bound. \end{abstract}

\maketitle

\section{Introduction}

A right-angled Artin group is a group given by a finite presentation whose only relations are commutators of the generators. These groups have nice algorithmic properties and act naturally on CAT(0) cube complexes.  Also known as graph groups, they occur in many different mathematical contexts; for some particularly interesting examples we refer to the work of Bestvina-Brady \cite{BesBra97} on finiteness properties of groups, Croke-Kleiner \cite{CroKle00} on boundaries of CAT(0) spaces, and Abrams \cite{Abr02} and Ghrist  \cite{Ghr01},\cite {GhrPet} on configuration spaces in robotics. For a general survey of right-angle Artin groups see \cite{Cha06}.

A nice way to describe a right-angled Artin group is by means of a finite simplicial graph $\G$.  If $V$ is the vertex set of $\G$, then the group $\AG$ is defined by the  presentation 
\[
\AG=\langle V \mid vw=wv \hbox{ if  $v$ and $w$ are connected by an edge in $\G$}\rangle.
\]
At the two extremes of this construction are the case of a graph with $n$ vertices and no edges, in which case $\AG$ is a free group of rank $n$, and that of a complete graph on $n$ vertices, in which case $\AG$ is a free abelian group of rank $n$.  In general, right-angled Artin groups can be thought of as interpolating between these two extremes. Thus it seems reasonable to consider automorphism groups of right-angled Artin groups as interpolating between  $Aut(F_n)$, the automorphism group of a free group, and $GL_n(\mathbb Z)$, the automorphism group of a free abelian group.  The automorphism groups of free groups and of free abelian groups have been extensively studied but, beyond work of Servatius \cite{Ser89} and Laurence \cite{Lau95}  on generating sets, there seems to be little known about the  automorphism groups of general right-angled Artin groups.  
 
 In this paper we begin a systematic study of automorphism groups of right-angled Artin groups.  We restrict our our attention to the case that the defining graph $\G$ is connected and triangle-free or, equivalently, $\AG$ is freely indecomposable and contains no abelian subgroup of rank greater than two.  A key example is when $\G$ is a  complete bipartite graph, which we call a {\it join\,} since $\G$  is the simplicial join of two disjoint sets of vertices.  The associated Artin group $\AG$ is a product of two free groups.  If  neither of the free groups is cyclic, then the automorphism group of $\AG$ is just the product of the automorphism groups of the two factors (or possibly index two in this product).  If  one of the free groups is cyclic, however, the automorphism group  is much larger, containing in addition an infinite group generated by transvections from the cyclic factor to the other factor. 
 
We will show that for any connected, triangle-free graph $\G$, the maximal join subgraphs of $\G$ play a key structural role in the automorphism group. We do this by proving that  the  subgroups $A_J$ generated by maximal joins $J$ are preserved up to conjugacy (and up to diagram symmetry) by automorphisms of $\AG$.  This gives rise to a homomorphism
\[
Out^0(\AG) \to \prod Out(A_J),
\]
where $Out^0(\AG)$ is a finite index normal subgroup of $Out(\AG)$ which avoids certain diagram symmetries.  We use algebraic arguments to prove that the kernel of this homomorphism is a finitely generated free abelian group. Elements of this kernel
commute with certain transvections, called ``leaf transvections," and adding them forms an even larger free abelian subgroup of $Out(\AG)$.  We show that this larger subgroup is the kernel of a homomorphism into a product of outer automorphism groups of free groups.  We then derive the Tits alternative for $Out(\AG)$ using the fact that the Tits alternative is known for outer automorphism groups of free groups.  

The second part of the paper takes a geometric turn.   The group $Out(F_n)$ can be usefully represented as symmetries of a topological space known as  {\it outer space}.  This space, introduced by Culler and Vogtmann in \cite{CulVog86},  may be described as a space of actions of $F_n$ on trees.  Outer space has played a key role in the study  of the groups $Out(F_n)$ (see for example, the survey article \cite{Vog02}).  For $GL(n,\Z)$ the analogous ``outer space" of actions of $\Z^n$ on $\R^n$ is the classical homogeneous space $SL(n,\R)/SO(n,\R)$. In Section~\ref{outerspace} of this paper we construct an outer space $\mathcal O(\AG)$ for the right-angled Artin group $\AG$ associated to any connected, triangle-free graph  $\G$.  If $\G$ is a single join,   $\mathcal O(\AG)$ consists of actions of $\AG$ on products of trees $T \times T'$.  In general, a point in outer space is a graph of such actions, parameterized by a collection of maximal joins in $\G$.  We prove that the space  $\mathcal O(\AG)$ is finite dimensional, contractible, and has a proper action of $Out^0(\AG)$. 

In the last section we give upper and lower bounds on the virtual cohomological dimension of $Out(\AG)$ and construct a spine for $\mathcal O(\AG)$, i.e. a simplicial equivariant deformation retract of $\mathcal O(\AG)$,  which realizes the most general upper bound.  A lower bound is given  by the rank of any free abelian subgroup, such as the subgroup found in the first part of the paper.  We show that this subgroup can be expanded even further to give a better lower bound.  
In  some examples, the upper and lower bounds agree, giving the precise virtual cohomological dimension of $Out(\AG)$.


\section{Preliminaries}
 

\subsection{Special subgroups}

A {\it simplicial graph} is a graph which is a simplicial complex, i.e.  a graph with no loops or multiple edges.  Vertices of valence one are called {\it leaves}, and all other vertices are {\it interior}.  To each finite simplicial graph $\G$  we  associate the {\it right-angled Artin group} $\AG$ as described in the introduction.  

A {\it special subgroup} of $\AG$ is a subgroup generated by a subset of the vertices of $\G$.   If  $\Theta$ is the full subgraph of $\G$ spanned by this subset, the special subgroup is naturally isomorphic to the Artin group $A_\Theta$.  In the discussion which follows, we will need to know the normalizers $N(A_\Theta)$, the centralizers $C(A_\Theta)$, and the centers $Z(A_\Theta)$ of  special subgroups $A_\Theta \subset \AG$.    To describe them, the following notation is useful.
\begin{definition} Let  $\Theta$ be a full subgraph of a simplicial graph $\G$. Then $\Theta\p$ is the intersection of the (closed) stars of all the vertices in $\Theta$:  $$\Theta\p=\bigcap_{v\in\Theta} st(v)$$
\end{definition}
Identifying the vertices of $\G$ with generators of $\AG$, one can also describe  $\Theta\p$ as the subgraph of $\G$ spanned by the vertices which commute with every vertex in $\Theta$.  We remark that this notation differs from that of Godelle \cite{God03}, who excludes points of $\Theta$ from $\Theta\p$.

\begin{proposition}\label{normalizer} For any graph $\G$, and any special subgroups
 $A_{\Theta}$ and $A_{\Lambda}$ of $\AG$,
\begin{enumerate}
\item the normalizer, centralizer, and center of $A_\Theta$ are given by
$$N(A_\Theta)=A_{\Theta\cup \Theta\p} \quad C(A_\Theta)=A_{\Theta\p}
\quad Z(A_\Theta)= A_{\Theta\cap \Theta\p}. $$
\item If $gA_{\Theta}g^{-1} \subseteq A_{\Lambda}$, then $\Theta\subseteq \Lambda$ and  $g=g_1g_2$ for some $g_1 \in N(A_{\Lambda})$, $g_2 \in N(A_\Theta).$
\end{enumerate}
\end{proposition}

\begin{proof} These statements are easily derived from work of Servatius  and work of Godelle, as follows. In \cite{Ser89}, Servatius proves that the centralizer of a single vertex $v$ is the special subgroup generated by $st(v)$ and hence for a set of vertices $\Theta$, the centralizer, $C(A_\Theta)$,  is generated by the intersection of these stars, which is exactly $\Theta\p$.  It follows that the center of $A_\Theta$ is $A_\Theta \cap C(A_\Theta)=A_{\Theta \cap \Theta\p}$.  
In \cite{God03}, Godelle considers normalizers and centralizers of special subgroups in a larger class of Artin groups,  Artin groups of ``FC type".  He defines the quasi-centralizer, $QZ(A_\Theta)$, of a special subgroup  to be the group of elements $g$ which conjugate the \emph{set} $\Theta$ to itself, and he proves that $N(A_\Theta)=A_\Theta \cdot QZ(A_\Theta)$.  In the right-angled case, no two generators are conjugate, hence $QZ(A_\Theta)=C(A_\Theta)=A_{\Theta\p}$ and $N(A_\Theta)=A_{\Theta\cup \Theta\p}$.

Godelle also describes the set of elements which conjugate one special subgroup $A_\Theta$ into another $A_\Lambda$ in terms of a category  $Ribb(V)$ whose objects are subsets of the generating set $V$ and whose morphisms conjugate one subset of $V$ into another.   In the case of a right-angled Artin group, since no two generators are conjugate, there are no morphisms between distinct objects of $Ribb(V)$ and the group of morphisms from an object $\Theta$ to itself is precisely the centralizer $C(A_\Theta)=A_{\Theta\p}$.  Proposition 3.2 of \cite{God03}  asserts that $gA_\Theta g^{-1} \subseteq A_\Lambda$ if and only if $g=g_1g_2$ where $g_1 \in A_\Lambda$ and $g_2$  is a morphism in a certain subcategory of $Ribb(V)$.  In the right-angled case,  it is straightforward to verify that such a morphism exists if and only if  $\Theta \subseteq \Lambda$ and $g_2 \in A_{\Theta\p}$.   (In Godelle's notation, he decomposes $\Theta$ into disjoint subsets $\Theta=\Theta_s \cup \Theta_{as}$ where $\Theta_s$ is the set of generators lying in the center of $A_\Theta$.  His theorem states that $g_2$ must commute with $\Theta_s$ and conjugate $\Theta_{as}$ to a set $R$ with $R \cup \Theta_s \subseteq \Lambda$.  In the right-angled case, this is possible only if  $R=\Theta_{as}$ and  $g_2$ also commutes with $\Theta_{as}$.)  Thus $g \in A_\Lambda A_{\Theta\p} \subseteq N(A_\Lambda)N(A_\Theta)$.
\end{proof}

\subsection{Cube complexes}
Associated to each right-angled Artin group $\AG$, there is a CAT(0) cube complex $C_\G$ on which $\AG$ acts, constructed as follows.  The 1-skeleton of $C_\G$ is the Cayley graph  of $\AG$ with generators the vertices $V$  of $\G$.   There is a cube of dimension $k>1$ glued in wherever possible, i.e. wherever the 1-skeleton of a cube exists in the Cayley graph.  In the quotient by the action there is a $k$-dimensional torus for each complete subgraph of $\G$ with $k$ vertices.
The cube complex associated to a free group is simply the Cayley graph of the free group, i.e. a tree on which the free group acts freely.  The cube complex associated to the complete graph on $n$ vertices is the standard cubulation of $ {\bf R}^n$.  The cube complex associated to a join $U\ast W$ is a product $T_U\times T_W$, where $T_U$ (resp $T_W$) is a tree on which the free group $F\langle U\rangle$ (resp $F\langle W\rangle$) acts freely with quotient a rose.  

\subsection{Generators for the automorphism group of a right-angled Artin group}

A set of generators for $Aut(\AG)$  was found by 
 M. Laurence \cite{Lau95}, extending work of H. Servatius \cite{Ser89}.  There are five classes of generators:
\begin{enumerate}
  \item  Inner automorphisms 
  \item  Inversions 
  \item   Partial conjugations   
   \item  Transvections
    \item  Symmetries  
\end{enumerate}

{\it Inversions\,} send a standard generator of $\AG$ to its inverse.

\begin{figure}[ht!]
\labellist
\small\hair 2pt  
\pinlabel {$v$} at 55 55
\pinlabel {$w_1$} [r] at -3 30
\pinlabel {$w_2$} [r] at -3 80
\endlabellist
\centering
\includegraphics[scale=0.6]{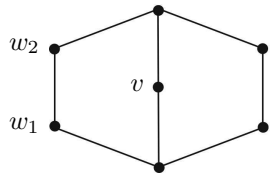}
\caption{Graph with a partial conjugation $w_i\mapsto v^{-1}w_iv$} 
\label{partconj}
\end{figure}

A {\it partial conjugation\,} exists when removal of the (closed) star of some vertex  $v$  disconnects the graph $\G$. In this case one obtains an automorphism by conjugating all of the generators in one of the components by $v$. (See example in Figure~\ref{partconj}.)

\begin{figure}[ht!]
\labellist
\small\hair 2pt 
\pinlabel {$v$} at 35 53
\pinlabel {$v$} at 307 42 
\pinlabel {$w$} at 149 46
\pinlabel {$w$} at 415 45
\pinlabel {Type I} [r] at 35 0
\pinlabel {Type II} [l] at   370 0
\endlabellist
\centering
\includegraphics[scale=0.8]{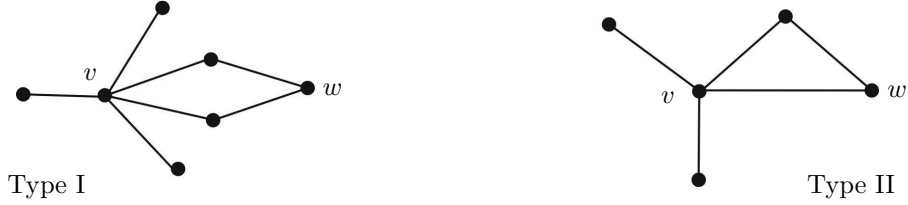}
\caption{Graphs with transvections $w\mapsto wv$} 
\label{transvections}
\end{figure}

{\it Transvections\,} occur whenever there are vertices  $v$ and $w$ such that $st(v)\supset lk(w)$; in this case the transvection sends $w\mapsto wv$.  There are two essentially different types of transvections, depending on whether or not $v$ and $w$ commute:
  \begin{enumerate}
  \item  {\it Type I transvections}.  $v$ and $w$  are not connected by an edge.
    \item {\it Type II transvections}.  $v$ and $w$  are connected by an edge.
\end{enumerate}
(See examples in Figure~\ref{transvections}.)

Finally, {\it symmetries\,} are induced by symmetries of the graph, and permute the generators.

We will be especially interested in the subgroup we obtain by leaving  out the graph symmetries:

\begin{definition} The  subgroup  of  $Aut(\AG)$   generated by inner automorphisms, inversions, partial conjugations and transvections is called the {\it pure automorphism group\,} and is denoted $Aut^0(\AG)$.  The image of $Aut^0(\AG)$ in $Out(\AG)$ is the group of {\it pure outer automorphisms\,} and is denoted $Out^0(\AG)$.
\end{definition}

The subgroups $Aut^0(\AG)$ and $Out^0(\AG)$ are easily seen to be normal and of finite index in $Aut(\AG)$ and $Out(\AG)$ respectively.  We remark that if $\AG$ is a free group or free abelian group, then $Aut^0(\AG)=Aut(\AG)$.  

\section{Maximal Joins}


\subsection{Restriction to connected, two-dimensional right-angled Artin groups}

If $\G$ is disconnected, then $\AG$ is a free product of the groups associated to the components of $\G$.  Guirardel and Levitt \cite{GuiLev05} have constructed a type of   outer space for a free product  with at least one non-cyclic factor, which can be used to reduce the problem of understanding the outer automorphism group to understanding the outer automorphism groups of the free factors.  Therefore,  in this paper we will consider only connected graphs $\G$. 

We will further restrict ourselves to the case that $\G$ has no triangles.  In this case,  the associated cube complex $C_\G$ is $2$-dimensional so we call these {\it two-dimensional} right-angled Artin groups.  To avoid technicalities, we also assume that $\G$ has at least 2 edges. 

Type II transvections are severely limited in the two-dimensional case.     Since $v$ and $w$ are connected by an edge, the vertex $w$ must actually be a terminal vertex i.e. a leaf: if there were another vertex $u\neq v$ connected to $w$, then the conditon $st(v)\supset lk(w)$ would imply that $u, v$ and $w$ form a triangle in the graph (see Figure~\ref{transvections}).  For this reason, we call Type II transvections {\it leaf transvections}.

If $\G$ is triangle-free and $\Theta$ is a subgraph with at least  one edge, then $\Theta\p\subset \Theta$, so by Proposition~\ref{normalizer},  
$$N(A_\Theta)=A_{\Theta} \hbox{  and } C(A_\Theta)=Z(A_\Theta)= A_{\Theta\p}.$$  If $\Theta$ contains two non-adjacent edges, then the latter groups are trivial.

\medskip
\noindent {\bf Key example.}  A join $\G=U\ast W$, has no triangles.    As we noted above, the associated right-angled Artin group $A_\G$ is $F\langle U \rangle\times F\langle W \rangle$.  It  is easy to deduce the structure of the automorphism group from Laurence's generators.   If $U$ and $W$ each contain at least two elements  then every automorphism preserves the two factors (or possibly switches them if $|U|=|W|$). Thus $Out(\AG)$ contains $Out(F\langle U \rangle) \times Out(F\langle W \rangle)$ as a subgroup of index at most 2.  If $U=\{u\}$ and $|W|=\ell\geq 2$, then $A_\G=\Z \times F\langle W \rangle$ with the center $\Z$ generated by $u$, and the elements of $W$ are all leaves. Any automorphism of $\AG$ must preserve the center and hence induces an automorphism of $F\langle W \rangle$, as well as an automorphism of the center $\Z$.  The  map $Out(\AG) \to Out(\Z)\times Out(F\langle W \rangle)$ splits and its kernel is the group generated by leaf transvections.  The leaf transvections commute, so $Out(\AG) \cong \Z^\ell  \rtimes (\Z/2 \times Out(F\langle W \rangle))$.  
\medskip

{\it We assume for the rest of this paper that that $\G$ is connected and triangle-free.  In addition, we assume that $\G$ contains at least two edges. }

\subsection{Restricting automorphisms to joins}

A connected, triangle-free graph $\G$ can be covered by subgraphs which are joins.  For example, to each interior vertex $v$ we can associate the join $J_v=L_v*L_v\p$, where
$L_v = lk(v)$. 
 Note that $L_v\p$ always contains $v$.  If $L_v\p=\{v\}$, then $J_v=st(v)$; in this case we say that $v$ is a {\it cyclic vertex} since $F(L_v\p)\cong \Z$ is cyclic.
 
 \begin{lemma} \label{edge}
If $v$ and $w$ are interior vertices joined by an edge of $\G$, then $L\p_v\subseteq L_w$, so we have

$$
\begin{array}{rccc}
   J_v= &L_v & \ast& L\p_v   \\
    &\cup  && \cap\\
    J_w=&L\p_w&\ast&L_w\\   
\end{array}
$$
In particular, $J_w \cap J_v = L\p_w \ast L\p_v$.
\end{lemma}

We remark that the $J_v$ is not properly contained in any other join subgraph of $\G$, i.e. $J_v$ is a {\it maximal join} in $\G$.  
The following proposition shows that the special subgroups $A_J$ associated to maximal join subgraphs $J$ of $\G$ are preserved up to conjugacy by pure automorphisms:  

\begin{proposition}\label{joins}
Let $\phi\in Aut^0(\AG)$ be a pure automorphism of $\AG$ and let $J=U \ast W$ be a maximal join in $\G$.  Then $\phi$ maps $A_J=F\langle U \rangle \times F\langle W \rangle$ to a conjugate of itself.  Moreover, if $U$ contains no leaves, then $\phi$ preserves the factor $F\langle U \rangle$  up to conjugacy. 
\end{proposition}

\begin{proof} It suffices to verify the proposition for the generators of $Aut^0(\AG)$.

{\it Inner automorphisms}.   These obviously send each $A_J$ to a conjugate of itself.

{\it Inversions}.  An inversion sends each $A_J$ to itself and preserves the factors. 

{\it Partial conjugations}. If  $\phi$ is a partial conjugation by a vertex $v$, we claim that $\phi$ either fixes all of $A_J$ or conjugates all of $A_J$ by $v$.  
Suppose first that $v$ is not in $U\ast W$.  The link of $v$ cannot contain vertices of both $U$ and $W,$ since there would then be a triangle in $\G$.  Furthermore, $lk(v)$ cannot contain all of $U$ since then adding it to $W$ would make a larger join, contradicting maximality.   Therefore the subgraph of $U\ast W$ spanned by vertices not in $lk(v)$ is still connected, so $\phi$  has the same effect (either trivial or conjugation by $v$) on generators corresponding to all vertices in  $U\ast W.$
Next, suppose that $v$ is actually in $U\ast W$, say $v\in U.$  The resulting partial conjugation restricted to  $A_{U\ast W}$ is an internal automorphism of  $A_{U\ast W}$ which may conjugate some generators of $\Fone$ by $v$, but has no effect on generators of $\Ftwo.$

{\it Transvections.}  We claim that a transvection either fixes $A_J$ or acts as an internal automorphism of $A_J$.  If $\phi$ is a transvection sending $s\to sv$, then $\phi$ is the identity on  $A_{U\ast W}$ unless $s\in U\ast W$, say $s\in U$.  If $s$ is not a leaf, then the condition $st(v)\supset  lk(s)$ and maximality imply that $v$ is also in $U$, so the restriction of $\phi$ is an internal automorphism of $A_{U\ast W}$ preserving the factor $\Fone$ and fixing $\Ftwo$.  If $s$ is a leaf, then $W=\{v\}$ and  $A_{U\ast W}= \Fone \times \Z$. In this case, $\phi$ fixes the (central) $\Z$ factor and multiplies $s$ by the generator of $\Z$.
\end{proof}

This proposition has two easy corollaries.  First, let 
$Sym(\G)$ denote the group of diagram symmetries of $\G$ and $Sym^0(\G)=Sym(\G) \cap Aut^0(\G)$.  Clearly $Aut(\AG)/Aut^0(\AG) \cong Out(\AG)/Out^0(\AG) \cong Sym(\G)/Sym^0(\G)$. Denote this quotient group by $Q(\G)$.  

\begin{corollary}\label{graphsymmetries}
The quotient maps from $Aut(\AG)$, $Out(\AG)$, and $Sym(\G)$ to $Q(\G)$ split.
\end{corollary}

\begin{proof}  It suffices to define a splitting of the projection $Sym(\G) \to Q(\G)$.  Composing with the inclusion of $Sym(\G)$ into $Aut(\AG)$ or $Out(\AG)$, gives a splitting in the other two cases. 

We first characterize elements of $Sym^0(\G)$ as those graph symmetries which only permute vertices with the same link.   
Define an equivalence relation on the vertices of $\G$ by $v \sim w$ if $lk(v)=lk(w)$.  
The elements of an equivalence class $[v]$ generate a free subgroup, and any automorphism of this subgroup extends (via the identity) to an automorphism of the 
whole Artin group $\AG$.  (These are the automorphisms generated by inversions and  transvections involving only elements of $[v]$.) Since $Aut^0=Aut$ for a free group, any permutation of $[v]$ can be realized by an element of $Sym^0(\G)$ which is the identity outside $[v]$.  Composing these gives an automorphism in $Sym^0(\G)$ realizing any permutation of the elements of each equivalence class.  

Conversely, if a graph symmetry $\phi$ is in $Sym^0(\G)$ we claim that it acts by permuting the elements of each equivalence class $[v]$. To see this, note that for any non-leaf vertex $v$, it follows from Proposition \ref{joins} that $\phi$ preserves $J_v$ and $L_v^\perp$, and hence it must also preserve $L_v$.  Since any graph symmetry takes links to links, $\phi$ permutes $[v]$.  Moreover, if $L_v$ contains a leaf $w$, then $\phi$ permutes $[w]$, the set of all leaves in $L_v$.

Now choose an ordering on the vertices in each equivalence class.  Then we can define a splitting of $Sym(\G) \to Q(\G)$  by mapping a coset to the unique element of the coset which is order preserving on every equivalence class.
\end{proof}

The second corollary of Proposition~\ref{joins} will be crucial for our analysis.

\begin{corollary}\label{restrict}
For every maximal join $J \subset \G$, there is a restriction homomorphism $R_J : Out^0(\AG) \to Out(A_J)$.
\end{corollary}

\begin{proof}  Fix $J$. Then for any element of $Out^0(\AG)$, there is a representative $\phi \in Aut^0(\AG)$ which maps $A_J$ to itself.  Any two such representatives differ by conjugation by an element of the normalizer of $A_J$.  But the normalizer $N(A_J)$ is equal to $A_J$, so the restriction of $\phi$ to $A_J$ is a well-defined element of $Out(A_J)$.  
\end{proof}

Let $\mathcal M$ be the set of all maximal join subgraphs of $\G$. We can put all of the homomorphisms $R_J$ for $J\in M$ together to obtain a homomorphism 
$$R=\prod_{J\in\mathcal M} R_{J} : Out^0(\AG) \to \prod_{J\in \mathcal M} Out (A_{J}).$$
To understand $Out^0(\AG)$, then, we would like to understand the image and kernel of this homomorphism.   But, first we note that there is a lot of redundant information in the set of all maximal joins used to define $R$; for example, the maximal joins of the form $J_v$ already cover $\G$.  

Even the covering of $\G$ by the maximal joins $J_v$ is inefficient.  If $v$ and $w$ have the same link, then $J_v=J_w$ so we don't need them both.  With this in mind, we now specify a subgraph $\G_0$ of $\G$ which will turn out to contain all of the information we need.  The key idea is that of vertex equivalence, which we already encountered in the proof of Corollary~\ref{graphsymmetries}.

\begin{definition} Vertices $v$ and $w$  of $\G$ are called {\it equivalent}  if they have the same link, i.e. $L_v=L_w$.  Equivalence classes of vertices are partially ordered by the relation $[v] \leq [w]$ if $L_v \subseteq L_w$. 
\end{definition}

To define $\G_0$, we choose a vertex in  each maximal equivalence class and let $\G_0$ be the full subgraph of $\G$ spanned by these vertices. 
In the special case that $\G$ is a star $\{v\}*W$, set $\G_0=\{v\}$.    
Up to isomorphism, $\G_0$ is independent of the choice of representatives. We denote by $V_0$ the set of vertices in $\G_0$.

\begin{figure}[ht!]
\labellist
\small\hair 2pt  
\pinlabel {$v$} at 2 22 
\pinlabel {$w$} at 68 22 
\endlabellist
\centering
\includegraphics[scale=1.2]{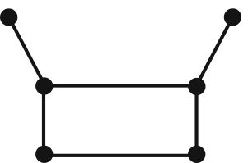}
\caption{$\G_0=[v,w]$}
\label{G0} 
\end{figure}

\begin{examples} 
(i)   If $\G$ is a tree then $\G_0$ is the subtree spanned by the vertices which are not leaves.

(ii)  If $\G$ is the graph in Figure~\ref{G0}, then $\G_0$ is the single edge spanned by $v$ and $w$.
\end{examples}

\emph{ Whether we are working with vertices of $\G$ or $\G_0$, the notation $L_v$  will always refer to the link  of $v$ in the original graph $\G$. Likewise,  a vertex is considered a leaf if it is a leaf of the full graph $\G$.  }

Recall that a vertex  $v$ of $\G$ is called a {\it cyclic vertex} if the associated maximal join $J_v$ is equal to  $st(v)$.  The following lemma specifies the properties of $\G_0$ which we will need.

\begin {lemma} \label{subgraph}
Let $\G_0$ be defined as above.  
\begin{enumerate}
\item $\G_0$ is a connected subgraph of $\G$.
\item The vertex set $V_0$ of $\G_0$ contains every cyclic vertex and no leaves of $\G$.  
\item Every vertex of $\G$ lies in $L_w$ for at least one $w \in V_0$ and lies in $L\p_w$ for at most one $w \in V_0$. 
\end{enumerate}
\end{lemma}

\begin{proof} (1) Let $v,w$ be two vertices in $\G_0$ and let
\[
v=v_0, v_1, \dots v_k=w
\]
be an edgepath in $\G$ connecting $v$ to $w$.  If $v_1$ does not lie in $\G_0$,  then there is a vertex $v_1' \in \G_0$ with $L_{v_1} \subseteq L_{v_1'}$.  Replacing $v_1$ by $v_1'$ gives another  edgepath in $\G$ from $v$ to $w$ whose first edge lies in $\G_0$.  The first statement of the lemma now follows by induction on $k$.

(2) Note that  $[v] \subseteq L\p_v$ for any $v$, and $[v]=L\p_v$ if and only if  $[v]$ is maximal.  Equality holds when $v$ is cyclic, since in this case $L\p_v=\{v\}$.  
If $v$ is a leaf then $v$ is connected to some  interior vertex $w$ by an edge (recall that we have assumed that the diameter of $\G$ is at least two). If $\G=st(w)$ we have defined $\G_0=\{w\}$ so $v\notin V_0.$  If $\G$ is not a star, $w$ must be connected to some other interior vertex $u$.  It follows that $L_v =\{w\} \subsetneq L_u$, hence $v \notin V_0$.

(3) Since $\G$ is connected, every vertex $v$ lies in the link of some other vertex and hence lies in the link of some maximal vertex.  Since $L\p_w=[w]$ for every vertex $w$ in $\G_0$, $v$ lies in at most one such $L\p_w$.
\end{proof}

\subsection{The kernel of the restriction and projection homomorphisms}\label{sec:kernel}

We are interested in determining the kernel of the map $R$ constructed from the restriction homomorphisms $R_J\colon Out^0(\AG)\to Out(A_J)$.  We first consider the kernel of the analogous map $R_0$ defined by looking only at the $R_J$ for maximal joins  $J=J_v$, for $v\in V_0$:    
\[
R_0=\prod_{v \in V_0} R_{J_v} : Out^0(\AG) \to \prod_{v \in V_0} Out (A_{J_v})
\]

This kernel consists of outer automorphisms such that any representative in  $Aut^0(\AG)$ acts by conjugation on each $A_{J_v},$ for $v\in \G_0$.

For any $\phi \in Out^0(\AG)$ and any interior vertex $v$, we can choose a representative automorphism $\phi_v$ such that $\phi_v (A_{J_v})=A_{J_v}$.  Since $L\p_v$ contains no leaves, $\phi_v$ must also preserve $A_{L\p_v}$ by Propositon~\ref{joins}. If $v$ and $w$ are interior vertices connected by an edge, then the representatives $\phi_v$ and $\phi_w$ are related as follows.

\begin{lemma}\label{edge-reps}  Suppose $v, w$ are interior vertices connected by an edge in $\G$, and $\phi\in Out^0(\AG)$ is represented by automorphisms  $\phi_v$ and $\phi_w$ with $\phi_v (A_{J_v})=A_{J_v}$ and $\phi_w(A_{J_w})=A_{J_w}$.
Then there exists $g_v \in A_{J_v}$ and $g_w \in A_{J_w}$ such that $c(g_v)\circ \phi_v = c(g_w) \circ \phi_w$, where   $c(g)$ denotes conjugation by $g$. 
\end{lemma}

\begin{proof}  Since $\phi_v$ and $\phi_w$ represent the same element in $Out^0(\AG)$, $\phi_w \circ \phi_v ^{-1}=c(g)$ for some $g \in \AG$.  Since $\phi_v$ preserves $A_{L\p_v}$, and $A_{L\p_v} \subset A_{J_w}$, we  have $g A_{L\p_v} g^{-1} \subset A_{J_w}$. By Proposition~\ref{normalizer} (2), we must then have $g=g_1g_2$ with $g_1 \in N(A_{J_w})=A_{J_w}$ and $g_2 \in N(A_{L\p_v})=A_{J_v}$.  Taking $g_w=g_1^{-1}$ and $g_v=g_2$ gives the desired formula.
\end{proof}

\begin{figure}[ht!]
\labellist
\small\hair 2pt
\pinlabel $v$ at 110 60
\pinlabel $w$ at  65 60
\pinlabel $u$ at 25 100
\endlabellist
\centering
\includegraphics[scale=0.7]{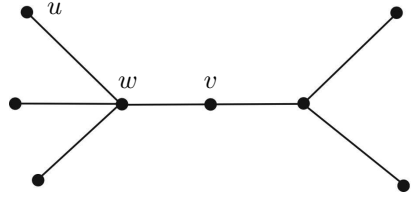}
\caption{$\G-st(v)$ has more components than $\G-\{v\}$} 
\label{separating}
\end{figure}

A vertex $v$ is a \emph{separating vertex} if $\G- \{v\}$ is disconnected.  It is easy to see that   separating vertices are cyclic.  It is also easy to see that conjugating any component of $\G- \{v\}$  gives an element of the kernel of $R$. We remark, however, that a component of $\G-\{v\}$ may contain more than one component of $\G-st(v)$, so not every partial conjugation by $v$ lies in this kernel. For example, in the graph in Figure~\ref{separating},  $\G-\{v\}$ has two components while $\G-st(v)$ has five, and the partial conjugation of $u$ by $v$ restricts non-trivially in $Out(A_{J_w})$.

In the rest of the paper, we will need to count several things associated to a vertex of $\G$, so we now establish some notation.

\begin{notation} Let  $v$ be a vertex of $\G$ and $\G_0$ the graph defined above. 
\begin{enumerate}
\item The valence, or degree, of $v$ in $\Gamma$ is denoted $\dv$.
\item If $v\in\Gamma_0$, the valence of $v$ in $\Gamma_0$ is denoted $\dov$.
\item The number of connected components of the complement $\Gamma-\{v\}$ is denoted $\cdv$.
\item The number of leaves attached by an edge to $v$ is denoted $\lv$.

\end{enumerate}
\end{notation}

\begin{proposition}\label{kernel}
The kernel $\kerR$ of the homomorphism 
$R_0$
is a finitely generated, free abelian group, generated by conjugations by  separating vertices $v$ on non-leaf components of $\G- \{v\}$.  If $\G$ is a star, the kernel is trivial; otherwise it has rank  $\sum_{v \in V_0} ( \cdv - \ell(v)-1)$.
\end{proposition}

\begin{proof}  If   $\G$ is the star of $v$ then $J_v=\G$  and $R_0$ is injective.  
So assume that $\G$ is not a star.  Note that in this case, $\AG$ has trivial center.

Let $r(v)=\cdv-\ell(v)$.  We will prove the theorem by defining a homomorphism  $$\mu=\prod_{v\in V_0} \mu_v\colon K_0\to \prod_{v\in V_0}\Z^{r(v)}/\Delta\iso  \prod_{v\in V_0}\Z^{r(v)-1},$$
where $\Delta$ is the diagonal subroup of $\Z^{r(v)}$, and showing that $\mu$ is an isomorphism.

Let $\phi$ be an element of  $K_0$.  For each $u\in V_0$  choose a  representative $\phi_u\in Aut^0(\AG)$ which acts as the identity on vertices of $J_u$. If $v$ and $w$ in $\G_0$ are connected by an edge, then $\phi_v^{-1}\phi_w$ is the identity on vertices of $ J_v\cap  J_w$.  Now $\phi_v^{-1}\phi_w$ is inner, so it is conjugation by a  (unique) element $g_{v,w}$ of the centralizer $C(A_{J_v\cap J_w})$.  By Proposition~\ref{normalizer} this centralizer is the special subgroup associated to $(J_v\cap J_w)\p$.  Since $\G$ has no triangles, $(J_v\cap J_w)\p=(L_v\p*L_w\p)\p$ is either the edge spanned by $v$ and $w$ (if both $v$ and $w$ are cyclic), the vertex $v$ (if $v$ is cyclic but $w$ is not), the vetex $w$ (if $w$ is cyclic but $v$ is not), or empty (if neither $v$ nor $w$ is cyclic).  Thus  $g_{v,w}=v^{-n}w^m$ with $n$ (resp. $m$) equal zero if $v$ (resp. $w$) is not cyclic, and   the equation of Lemma~\ref{edge-reps} reduces to
\[
c(v^{n})\circ \phi_v=c(w^m)\circ \phi_w.
\]

We are now ready to  define $\mu_v\colon K_0\to \Z^{r(v)}/\Delta$.  Let $C_1,\ldots,C_{r(v)}$ be the non-leaf components of  $\G-\{v\}$.  In each $C_i$, choose a vertex $w_i$ adjacent to $v$ which is also in $\G_0$.  We have $c(v^{n_i})\circ \phi_v=c(w_i^{m_i})\circ \phi_{w_i}$ for unique integers $n_i$ and $m_i$ and we define  $\mu_v(\phi)=[(n_1,\ldots,n_{r(v)})]$. 

We have made several choices, and we must show that $\mu_v(\phi)$ is independent of these choices.  The integers $n_i$ depend, a priori,  on our choice of representatives $\phi_v$ and $\phi_{w_i}$.  But $\phi_v$ is unique up to conjugation by an element of $\AG$ centralizing $J_v$, namely a power of $v$, and similarly for $\phi_{w_i}$.  It follows that the choice of $\phi_{w_i}$ has no effect on $n_i$ whereas the choice of a different representative for $\phi_v$ will change all of the integers $n_i$ by the same amount.  Therefore  the class of $(n_1,\ldots, n_{r(v)})$ modulo the diagonal $\Delta$ is independent of $\phi_v$ and the $\phi_{w_i}$.  

We now show that $\mu_v(\phi)$ is independent of the choice of the $w_i$.  Let $w_i'$ be a different choice, and connect $w_i$ to $w_i^\prime$ by a simple path in $C_i$, i.e. a path which goes through each vertex at most once. In particular, each vertex on the path is an interior vertex.    

 We first claim that the projection of $\phi$ in $Out(J_u)$ is trivial for every interior vertex of $\G$, not just for those in $\G_0$.   Suppose $u$ is an interior vertex of $\G-\G_0$.  By Lemma~\ref{subgraph},  $u$ is adjacent to some $w \in \G_0$, so $L_u^\perp \subset L_w$. By definition of $\G_0$,  $L_u \subset L_v$ for some $v \in \G_0$, so $J_u \subset L_v \cup L_w$.   Since $v$ and $w$ are in $\G_0$ and connected by an edge, we can find $n$ and $m$ with $c(v^{n})\circ \phi_v=c(w^m)\circ \phi_w$.  Setting $\phi_u=c(v^{n})\circ \phi_v$ we have that $\phi_u$ acts trivially on every vertex of $J_u$.   

  For any two interior vertices $u,u'$ of $\G$, we define $g_{u,u'}$ to be the (unique) element of $\AG$ such that $\phi_u\circ \phi_{u'}^{-1}$ is conjugation by $g_{u,u'}$.  If $u$ and $u'$ are connected by an edge, then in light of the previous paragraph, the same argument used for the case $v,w\in\G_0$ applies to show that $g_{u,u'}=u^{-n}(u')^m$ for some integers $m$ and $n$.  
  
  We now return to the simple path in $C_i$ joining $w_i$ and $w_i^\prime$.  For each edge $[u,u']$ of this path,   $g_{u,u'}$ is a word which does not involve $v$.  Observe that $g_{w,w'}$ is the product of the $g_{u,u'}$, so that  $g_{w,w'}$ does not involve $v$.  It follows that the powers of $v$ in  $g_{v,w_i}$ and $g_{v,w_i'}$  are the same, showing that $\mu(\phi)$ is well-defined.  
  
It is straightforward to verify that $\mu_v$ is a homomorphism.

To see that $\mu_v$ is surjective, take any $r(v)$-tuple of integers $(n_1,\ldots, n_r)$.  The product  $\phi$ of partial conjugations of $C_i$ by $v_i^{n_i}$ satisfies $\mu_v(\phi)=(n_1,\ldots, n_r)$. 

Finally, we show that $\mu=\prod{v\in V_0} \mu_v$ is injective.  
Suppose $\phi$ lies in the kernel of $\mu$.  Then we can choose a representative $\phi_v$ such that for any $w$ adjacent to $v$,  $g_{v,w}$ is just a power of $w$.  The same reasoning applied to $w$ implies that $g_{w,v}$ is just a power of $v$.  But $g_{v,w}=g_{w,v}^{-1}$ so we must have $g_{v,w}=1$. It  follows that for any adjacent pair of vertices, $\phi_v=\phi_w$.  Since $\G$ is connected, this gives a representative of $\phi$ which acts trivially on the join of every vertex; in other words, $\phi$ is trivial in $Out^0(\AG)$.
\end{proof}

\begin{remark}  Since the generators of $K_0$ given by Proposition~\ref{kernel} restrict to inner automorphisms on \emph{every} join, it follows from the theorem that   the homomorphism $R$, which was defined over all maximal joins $J$ instead of just the joins $J_v$ with $v \in V_0$, has the same kernel as $R_0$. 
\end{remark}

One advantage of restricting attention to the joins $J_v$ for $v\in V_0$ is that we can further define a {\it projection homomorphism}, as follows.
Since vertices of $V_0$ are interior, $L_v\p$ contains no leaves, hence by Proposition \ref{joins}, every $\phi \in Out^0(\AG)$  has a representative $\phi_v$ which preserves both $A_{J_v}$ and $A_{L_v^\perp}$.  Thus $\phi_v$ descends to an automorphism $\bar\phi_v$ of $A_{L_v}=A_{J_v} / A_{L_v^\perp}$.  This gives rise to a homomorphism  $P_v : Out^0(\AG) \to Out(F\langle L_v\rangle)$.  Let $P$ be the product homomorphism
\[
P = \prod_{v \in V_0} P_v : Out^0(\AG) \to \prod Out(F\langle L_v\rangle). 
\]

Recall that  $\cdv$ denotes the number of connected components of  $\Gamma-\{v\}$.

\begin{proposition}\label{KerP}   The kernel $\kerP$ of $P$ is a free abelian group, generated by $\kerR$ and the set of leaf transvections.  If $\G$ is a star $\{v\}\ast W$ then $K_P$ has rank $|W|$; otherwise, it has   rank $ \sum_{v \in V_0} (\cdv-1)$. \end{proposition}

\begin{proof}  Let $\ell$ be the number of leaves in $\G$.  It is clear that leaf transvections are contained in $\kerP$ and that they generate a free abelian group of rank $\ell$. It is also easy to see that leaf transvections commute with the generators of  $\kerR$ since if $u \mapsto uv$ is a leaf transvection, then $u$ and $v$ are connected by an edge and hence belong to the same component of $\G -\{w\}$ for any $w \in \G_0$. Together with $\kerR$, the leaf transvections thus generate a free abelian subgroup of the specified rank, by Proposition~\ref{kernel}.

It remains only to show that this subgroup is all of $\kerP$. 
Consider an element $\phi$ in $\kerP$. For $v \in V_0$, let $\phi_v$ denote a representative automorphism which preserves   $A_{J_v}$ and hence also $A_{L_v^\perp}$.  
Then for $w \in L_v$, we have $\phi_v(w)=wg$ for some $g \in A_{L_v^\perp}$ (up to conjugation by an element of $A_{J_v}$). Since $w$, and hence $\phi_v(w)$,  commutes with $A_{L_v^\perp}$, $g$ must lie in the center of $A_{L_v^\perp}$.  If $v$ is not cyclic, the center is trivial, so $g=1$.  If $v$ is cyclic, then then $g=v^k$.   Any automorphism preserve centralizers, so if $k \neq 0$,  $\phi_v(C(w))=C(wv^k)=\langle w, v \rangle$, so  $C(w)$ is free abelian of rank 2. This implies that $w$ is a leaf.  We conclude that $\phi_v$ acts as the identity on non-leaf elements of $L_v$ and as leaf transvections on leaf elements. It follows that there exists a product of leaf transvections $\theta$ such that $\phi \circ \theta$ lies in $K_0$. 
\end{proof}

We conclude this section with an easy consequence of Proposition~\ref{KerP}.  We say that a group $G$ satisfies the \emph{Tits alternative} is every subgroup of $G$ either contains a non-abelian free group or is virtually solvable. Tits \cite{Tit72} proved that all finitely generated linear groups satisfy the Tits alternative and Bestvina, Feighn, and Handel \cite{BFH00}, \cite{BFH05} proved that $Out(F_n)$ does likewise.

\begin{theorem}
If $\G$ is connected and triangle-free, then $Out(\AG)$ satisfies the Tits alternative.
\end{theorem}

\begin{proof}  It is an easy exercise to check that the property that a group satisfies the Tits alternative is preserved under direct products and abelian extensions.  Since $Out(F\langle L_v\rangle)$ satisfies the Tits alternative by Bestvina-Feign-Handel, it follows that the image of $P$ also does.  Since the kernel of $P$ is abelian, we conclude that $Out^0(\AG)$, and hence also $Out(\AG)$, satisfies the Tits alternative.  
\end{proof}

\begin{remark}
Many of the results in this section generalize to higher dimensional right-angled Artin groups.  The details will appear in a subsequent paper.
\end{remark}

\section{Outer Space}\label{outerspace}


In this section we introduce ``outer space" for a right-angled Artin group $\AG$. We continue to assume that $\G$ is a connected, triangle-free graph and has diameter $\geq 2$. 

Let $\mathcal O(F)$ denote the unreduced, unprojectivized version of Culler and Vogtmann's outer space for a free group $F$ \cite{CulVog86}.  This space can be described as the space of minimal, free, isometric actions of $F$ on simplicial trees. (The terms ``unreduced" and ``unprojectivized" specify that quotient graphs may have separating edges, and that we are not considering homothetic actions to be equal.)    Our initial approach to constructing outer space for $\AG$ was to consider minimal, free, isometric actions of $\AG$ on  CAT(0) 2-complexes.  In the case of a single join, these turn out to be products of trees.  More generally, (under mild hypotheses) such a 2-complex is a union of geodesic subspaces which are products of trees.  However, the interaction between these subspaces proved difficult to control and we ultimately found that it was easier to work directly with the tree-products.

\subsection{Outer space for a join}\label{OSjoin}
Let us examine more closely the case when $\G=U\ast W$ is a single join.  
Suppose that $\AG=F\langle U \rangle\times F\langle W \rangle$ acts freely and cocompactly by isometries on a piecewise Euclidean CAT(0) 2-complex $X$ with no proper invariant subspace.

If neither $U$ nor $W$ is a singleton, then $\AG=F\langle U \rangle\times F\langle W \rangle$ has trivial center, so the splitting theorem for CAT(0) spaces (\cite{BriHae99}, Theorem 6.21) says that $X$ splits as a product of two one-dimensional CAT(0) complexes  (i.e. trees) $T_U \times T_W$, and the action of $\AG$ is {\it orthogonal}, i.e. it is the product of the actions of   $F\langle U \rangle$  on $T_U$ and $F\langle W \rangle$ on $T_W$.  Twisting an action by an element of $Out(\AG)$, which contains $Out(F\langle U \rangle)\times Out(F\langle W \rangle)$ as a subgroup of index at most two, preserves the product structure. Thus it makes sense to take as our outer space the product of the Culler-Vogtmann outer spaces $\mathcal O(F\langle U \rangle) \times \mathcal O(F\langle U \rangle)$.

If $U=\{v\}$ is a single vertex   then $\AG=\Z\times F\langle W \rangle$ where $\Z$ is generated by $v$.  In this case, $X$ is equal to the min set for $v$, so that $X$  splits as a product $\alpha_v\times T_W$, where $\alpha_v$ is an axis for $v$ and $T_W$ is a tree  (\cite{BriHae99}, Theorem 6.8) .  We identify the axis $\alpha_v$ with a real affine line, with $v$ acting by translation in the positive direction by an amount $t_v>0$.  The tree $T_W$ has a free $F\langle W \rangle$-action induced by the projection $\AG\to F\langle W \rangle$, but the action of $\AG$ on $X$ need not be   an orthogonal action.  Recall that $Out(\AG) \cong \Z^\ell \rtimes (\Z/2 \times Out(F\langle W \rangle)$ where $\ell=|W|$ and $\Z^\ell$ is generated by leaf transvections.  Twisting an orthogonal action by a leaf transvection $w\mapsto wv$   results in an action which is no longer orthogonal.  Instead, $w$ now acts as translation in a diagonal direction on the plane in $\alpha_v\times T_W$ spanned by   $\alpha_v$ and the axis for $w$ in $T_W$, i.e. the translation vector has a non-trivial $\alpha_v$-component.   More generally, for any free minimal action of $\AG$ on a CAT(0) 2-complex $X=\alpha_v\times T_W$, the generator 
$v$ acts only in the $\alpha_v$-direction,
$$v\cdot (r,x) = (r+ t_v, x),$$
while an element $w_i \in W$  has a ``skewing constant" $\lambda(w_i)$, i.e. $w_i$  acts by 
$$w_i\cdot (r,x)=(r+\lambda(w_i), w_i\cdot x).$$
A free action of $\AG$ on $X$ is thus determined by an $F\langle W \rangle$-tree $T_W$, the translation length $t_v$ and an $\ell$-tuple of real numbers $(\lambda(w_1),\ldots,\lambda(w_\ell))$.   So it is reasonable in this case, to take for outer space the product $\mathcal O(F\langle W \rangle)\times \R_{>0}\times \R^\ell $.

\subsection{Tree spaces}

Our analysis of the join case motivates the following definition.

Let $J=U \ast W$ be a join in $\G$.  An \emph{admissible tree-space} $X_J$ for $J$ is a product of two simplicial, metric trees, $T_U$ and $T_W$ with free, minimal, isometric actions of $\Fone$ and $\Ftwo$ respectively, and an action of $A_J=\Fone \times \Ftwo$ on $X_J= T_U \times T_W$ of the following type.
\begin{enumerate}
\item  If  $J$ contains no leaves, then the action is the product of the given actions, 
\[
(g_1, g_2) \cdot (x_1, x_2) = (g_1\cdot  x_1, g_2\cdot  x_2).
\]
\item  Suppose $J$ contains leaf vertices, say in $W$.  This forces $U=\{v\}$, so  $A_J=\langle v\rangle\times F\langle W\rangle$, $T_U=\R$ and $v$ acts on $T_U$ as translation by some positive real number $t_v$.   Then there exists a homomorphism  
$\lambda : \Ftwo  \to \mathbb R$ which is zero on non-leaf vertices of $W$, such that
\[
 (v^n, g) \cdot (r , x) = (r+n t_v+\lambda (g), \, g \cdot   x).
 \]
\end{enumerate}
We remark that the definition of admissible depends not only on the join $J$, but on the graph $\G$ as well since $\G$ determines which vertices are considered as leaves.  Since $\G$ is fixed throughout, this should not cause any confusion.

A point in outer space for $\AG$ will be a collection of admissible tree-spaces satisfying certain compatibility conditions.  Recall that to each interior vertex of $\G$ we have associated a maximal join $J_v=L_v \ast L_v\p$.   If $e$ is an edge from $v$ to $w$, set $J_e=J_v \cap J_w= L\p_v \ast L\p_w$. If $e$ lies in $\G_0$, then $J_e$ contains no leaves.

\begin {definition}  A \emph{graph of tree-spaces} $\mathfrak X=\{X_v,X_e,i_{e,v}\}$ for $A_\G$ consists of the following data.
\begin{enumerate}
\item  For each vertex $v \in \G_0$, an admissible tree-space $X_v$ for $J_v$,
\item  For each edge $e \in \G_0$ with vertices $v$ and $w$, an admissible tree space $X_e$ for  $J_e$ and a pair of  $A_{J_e}$-equivariant isometric  embeddings.
\[
X_v \overset {i_{e,v}}  \longleftarrow     X_e  \overset {i_{e,w}}  \longrightarrow X_w.
\]
\end{enumerate}

We define \emph{outer space} for $A_\G$ to be the set
\[
\OG = \{ \mathfrak X \mid \textrm {$ \mathfrak X$ is a graph of tree-spaces for $A_\G$} \} / \sim
\]
where $\sim$ is the equivalence relation induced by replacing any $X_v$ (respectively $X_e$) with an equivariantly isometric space $X_v'$ (respectively $X_e'$), and composing the associated connecting maps $i_{e,v}$, by the equivariant isometry. A natural topology for $\OG$ will be described in Section~\ref{topology}.
\end{definition}

Since $J_e$ contains no leaves, the group $A_{J_e}=F\langle L\p_v\rangle \times F\langle L\p_w\rangle$ acts orthogonally on 
$X_e$ and the maps $i_{e,v}$ and $i_{e,w}$ split as products. 
Write $X_v=T_v \times \TvP$, where $T_v$ is an $F\langle L_v\rangle$-tree, and $\TvP$ is an $F\langle L_v\p\rangle$-tree.  Then, up to equivariant isometry, we may assume that $X_e= \TvP \times \TwP$, and that $i_{e,v}=(i_1,i_2)$ is the identity on the first factor while $i_{e,w}=(j_1,j_2)$ is the identity on the second.  Denote the embedding $i_2: T_w^\perp \to T_v$ by $i(w,v)$.  The \emph{image} of $i(w,v)$ is uniquely determined, namely it is the minimal $F\langle L\p_w\rangle$-invariant subtree of $T_v$.  However, the \emph{map} $i(w,v)$  is not necessarily unique.  
It is unique if $w$ is not cyclic, since in this case $\TwP$ has no non-trivial equivariant isometries. But if  $w$ is cyclic  then $\TwP$ is a real line and the possible equivariant inclusions $i_2$ are parameterized by $\mathbb R$.

 \subsection{Basepoints} \label{Xbpts}  
 To keep track of the inclusions $i(w,v)$, it will be convenient to introduce basepoints.  These will also play a crucial role in the proof of the contractibility of $\OG$.

 {\it Basepoints for free actions of free groups.}
  Let $F\langle S\rangle$ be a finitely generated free group with a specified basis $S$ of cardinality at least 2, and let $T$ be a metric tree with a minimal, free, isometric action of $F\langle S\rangle$.  Each generator $s \in S$ preserves a unique line $\alpha(s)$ in $T$ called the axis for $s$.  Orient the axis so that $s$ acts as translation in the positive direction. We  choose a base point $b(s)$ on $\alpha(s)$ as follows.  
 
\begin{figure}[ht!]
\labellist
\small\hair 2pt  
\pinlabel {$\alpha(t)$} [r] at 115 112 
\pinlabel {$\alpha(s)$} [l] at 253 68
\pinlabel {$\alpha(t^\prime)$} [r] at 82 6 
\pinlabel {$\alpha(t^{\prime\prime})$} at 120 20 
\pinlabel {$p(t,s)$} [l] at 170 77 
\pinlabel {$p(t^\prime,s)$} [r] at 142 77 
\pinlabel {$p(t^{\prime\prime},s) = b(s)$} at 45 77 
\endlabellist
\centering
\includegraphics[scale=0.8]{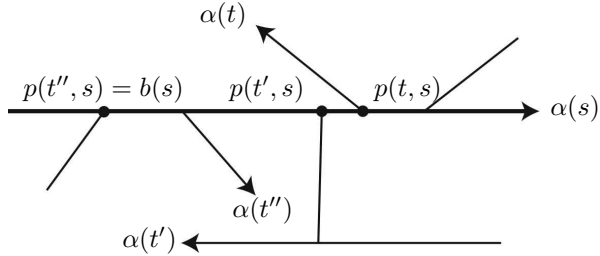}
\caption{Projections and basepoint on the axis   $\alpha(s)$} 
\label{projections}
\end{figure}

 For each generator $t \neq s$, the set of points on $\alpha(s)$ of minimal distance from $\alpha(t)$ is a closed connected interval (possibly a single point). Define the projection $p(t,s)$ of $\alpha(t)$ on $\alpha(s)$ to be the initial point of this interval, and the basepoint $b(s)$ to be the minimum of these projections, with 
  respect to the ordering given by the orientation of $\alpha(s)$.  (See Figure~\ref{projections}).  The basepoints $b(s)$ will be called the {\it unrestricted basepoints} of the action. 

\vskip 2pt  

{\it Basepoints in $\X$.}  Suppose now that we have a graph $\X$ of tree-spaces for $\AG$.  
Assume $\G$ is not a star, so $\G_0$ contains at least two vertices.  If $L_w\p$ contains at least two vertices, then for each $u\in L_w\p$ we let $b(u)$ be the unrestricted basepoint  on the axis $\alpha(u)$  for the action of $F\langle L_w\p\rangle$ on $T_w\p$.

  If $w$ is cyclic the axis $\alpha(w)$ is the  {\it entire} tree $\TwP$,  in which case we cannot use this method to choose a basepoint.    If $v$ is adjacent to $w$ in $\G_0$, we have an equivariant isometry $i(w,v)\colon \TwP\to  T_v$.  In $T_v$, the element $w$ has an axis $\alpha_v(w)$ with an unrestricted basepoint $b_v(w)$.  We can use $i(w,v)$ to pull this back to a point $c_v(w)$  on $\alpha(w)$.   
  
\begin{figure}[ht!]
\labellist
\small\hair 2pt
\pinlabel $w$ at 128 335
\pinlabel $\Gamma$ at 100 355
\pinlabel $v$ at 63 335
\pinlabel  $u$ [l] at 200 350 
\pinlabel  $u^\prime$ [l] at 200 300
\pinlabel  $T_{v}$ at 45 50
\pinlabel  $T_{w}^\perp$ at 225 50
\pinlabel  $T_{u}$ at 365 265
\pinlabel  $T_{{u'}}$ at 365 16
\pinlabel  $i(w,v)$ at 142 171
\pinlabel $i(w,u)$ [r] at 305 210
\pinlabel $i(w,u')$ [l] at 303 125
\pinlabel $b_v(w)$ at 24 146
\pinlabel $b_u(w)$ [r] at 417 210
\pinlabel $b_{u'}(w)$ [r]  at 417 87
\pinlabel $\alpha_v(w)$   at 45 247
\pinlabel $\alpha(w)$ at 225 245
\pinlabel $\alpha_u(w)$ at 420 355
\pinlabel $\alpha_{u'}(w)$  at 420 153
\pinlabel $c_v(w)$ [r] at 220 140 
\pinlabel $c_u(w)$ [l] at 230 168
\pinlabel $c_{u'}(w)=b(w)$ [l] at 230 99
\endlabellist
\centering
\includegraphics[scale=0.8]{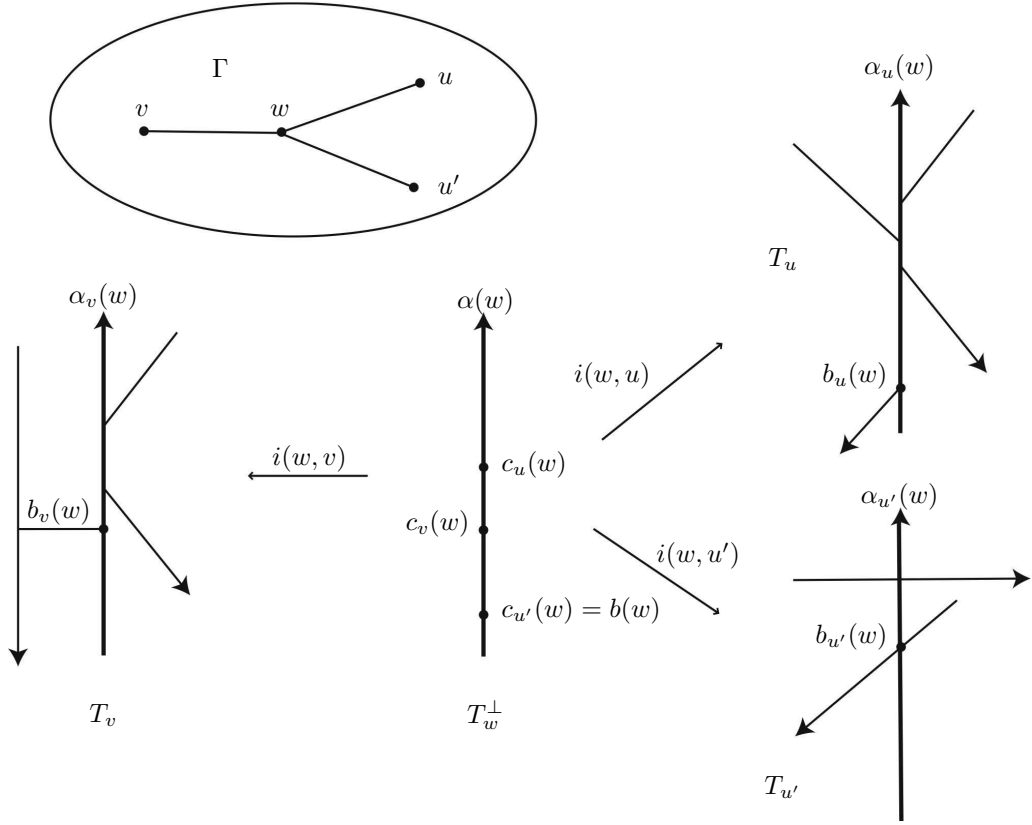}
\caption{Basepoint on the axis for $w$ in $T_w^\perp$} 
\label{axes}
\end{figure}

  In this way, we get a point $c_v(w)$ on $\alpha(w)$ for each $v$ adjacent to $w$ in $\G_0$.  We take the minimum $b(w)$  of these points as a basepoint for $\alpha(w)$. (See Figure~\ref{axes}.)
   
\subsection{Topology on $\OG$}\label{sec:topology}

We will define the topology on $\OG$ by embedding $\OG$ in a product of topological spaces  and giving it the resulting subspace topology.  

In the case that $\G$ is a star $\{v\}*W$,  we have seen that $\OG$ can be identified with $\mathcal O(F\langle W \rangle) \times \R_{>0} \times \R^\ell$, where $\ell=|W|$.  We take this to be a homeomorphism.

For any graph $\G$, recall that $V_{cyc}$ denotes the set of cyclic vertices and $\dov$ the valence of $v$ in $\G_0$.  

\begin{proposition}\label{topology}  Let $\G$ be a graph which is not a star. Let $\ell$ be the number of leaves of $\G$, and let  $k=\sum_{v \in V_{cyc}} \dov$.  Then there is an injective map 
\[
\OG\hookrightarrow \left(\prod_{v \in \G_0} \mathcal O(F\langle L_v\rangle) \right) \times  \mathbb R^k \times \mathbb R^\ell.
\]
The image of this map is of the form $Y \times Q \times \R^\ell$,  where $Y$ is a subspace of $\prod_{v \in V_0} \mathcal O(F\langle L_v\rangle)$ and $Q$ is a piecewise linear subspace of $ \R^{k}$ of dimension $\sum_{v \in V_{cyc}} (\dov-1)$.
\end{proposition}

\begin{proof}  
 Let $\X$ be an element of  $\OG$.  
The vertex spaces $X_v$ are determined by free group actions on two trees, $\TvP$ and $T_v$, together with a real number $\lambda(w)$ for each leaf vertex $w \in L_v$.  Since $\G$ is not a star and $\G_0$ is connected, there is at least one vertex $w$ of $L_v$ which is also in $\G_0$.  The tree  $\TvP$   is equivariantly isometric to the minimal $F\langle L\p_v\rangle$ subtree of $T_w$ for this $w$.  Thus the adjacency relations of $\G$, the trees $T_v$ for $v\in V_0$  and the real numbers $\lambda(w)$ for $w$ a leaf completely determine the vertex spaces $X_v$.  As noted above, the edge spaces are of the form $X_e = \TuP \times \TvP$ with orthogonal action of $F\langle L\p_u\rangle \times F\langle L\p_v\rangle$, so that they too are determined by the actions on the trees $T_v$.

It remains to account for the connecting maps between edge spaces and vertex spaces. If $w$ is not cyclic,  then for any vertex $v\in L_w\cap  \G_0$, there is a unique equivariant map from $\TwP$ into $T_v$, so the connecting map  $X_e \to X_v$ is uniquely determined.   If  $w$ is cyclic, then $\TwP$ is a copy of the real line, and an equivariant embedding of $\TwP \to T_v$ is determined by the position of the image of $b(w)$ on $\alpha_v(w)$ i.e. by the difference $\Delta_v(w)=c_v(w)-b(w)$, where the points $c_v(w)$ and $b(w)$ are as defined in Section~\ref{Xbpts} and illustrated in Figure~\ref{axes}. 

Associated to  $\X$ we have a point $T_v \in \mathcal O(F\langle L_v\rangle)$ for each $v \in V_0$, a real number $\lambda(w)$ for each leaf $w$, and a $\delta_0(w)$-tuple of real numbers $(\Delta_v(w))$, for each $w \in V_{cyc}$.  This data completely determines $\X$.  Thus, we have an injective map
\[
f: \mathcal O_\G \to \left(\prod_{v \in V_0} \mathcal O(F\langle L_v\rangle) \right) \times  \mathbb R^{k} \times \mathbb R^\ell.
\]

Now consider the image of $f$.  There are no restrictions at all on the numbers $\lambda(w)$; given any $\X$ in $\OG$, arbitrarily changing $\lambda(w)$ for any leaf $w$ gives rise to another valid graph of tree-spaces.  As for the numbers $\Delta_v(w)=c_{v}(w)-b(w)$, since $b(w)$ was defined to be the infimum of the $c_{v}(w)$, we must have $\Delta_{v}(w) \geq 0$ for all $v  \in L_w \cap V_0$, and at least one of these must equal $0$. There are no other restrictions;  the points $c_{v}(w)$ can be varied independently of each other by changing a single embedding $i(w,v)$.  

Let $Q_w \subset \mathbb R^{\delta_0(w)}$ be defined by 
\[
Q_w=\{(r_1, \dots r_{\delta_0(w)}) \mid \textrm{$r_i \geq 0$ for all $i$ and $r_j=0$ for some $j$}\}.
\]
In other words, $Q_w$ is the boundary of the positive orthant of $\R^{\delta_0(w)}$.
Then the image of $f$ is of the form $Y  \times \prod_{w \in  V_{cyc}} Q_w\times \R^\ell$ where $Y$ is a subspace of $\prod_{v \in V_0} \mathcal O(F\langle L_v\rangle)$.  
\end{proof}

\begin{remark}\label{Q}
For future reference, we remark that $Q_w$ can also be identified with $\R^{\delta_0(w)}/\R$ (where $\R$ acts as diagonal translation) with coordinates given by the $c_{v}(w)$'s.  
\end{remark}

\begin{example}\label{join}  
If $\G = U\ast W$ is a join with $|U|,|W| \neq 1$, then $\G_0$ consists of a single edge $e$ joining a pair of vertices $u \in U$ and $w \in W$.  In this case, $J_e=J_u=J_w=\G$ and  a point in $\mathcal O(\AG)$ is determined by a single tree space for $\G$ (since there is a unique equivariant isometry of such a tree-space).  In this case, $\ell=\delta_0(u)-1=\delta_0(w)-1=0$ and  $Y$ is all of $\mathcal O(F\langle U \rangle) \times \mathcal O(F\langle W \rangle)$.   Thus, the definition of outer space for a join agrees with that proposed  in Section~\ref{OSjoin}. 
\end{example}

\emph{From now on,  $\OG$ will be viewed as a topological space with the subspace topology induced by the embedding $f$ described above.}  Each space $\mathcal O(F)$ is endowed with the equivariant Gromov-Hausdorff topology.  In this topology, a neighborhood basis of an action of $F$ on a tree $T$ is given by the sets $N(X,H,\epsilon)$, where $X=\{x_1,\ldots,x_N\}$ is a finite set of points in $T$ and $H=\{g_1,\ldots,g_M\}$ is a finite set of elements of $F$.  An action of $F$ on $T'$ is in this neighborhood if there is a subset $X'=\{x_1',\ldots,x_N'\}$ of $T'$ such that 
$$|d_T(x_i,g_jx_k)-d_{T'}(x'_i,g_jx'_k)|<\epsilon$$
for all $i,j,k.$


\section{Contractibility}

In this section we prove the following theorem:

\begin{theorem}\label{contractible}
  For any connected, triangle-free graph $\G$, the space $\OG$ is contractible. 
\end{theorem}

The proof follows ideas of Skora \cite{Sko} and Guirardel-Levitt \cite{GuiLev05, GuiLev06} on  ``unfolding" trees.  In their work, however, the unfolding of a tree was defined with respect to a single basepoint.  In our case we will need to preserve several basepoints.  

We first consider a single action.  
Let $F\langle S\rangle$ be a free group with a preferred generating set $S$ and let $T \in \mathcal O(F\langle S\rangle)$ be an $F\langle S\rangle$-tree.

\begin{proposition}\label{BT} Let $P=\{p(t,s)\}$ be any set of projections, with  $s,t\in S$. The convex hull $B(P,T)$ of $P$ in $T$ depends continuously on $T$.
\end{proposition}
\begin{proof}
If $T'$ is any other $F\langle S\rangle$-tree, let $P'=\{p'(t,s)\}$ be the corresponding projections in $T'$.  We must show:   given $\epsilon>0$ and $T\in \mathcal O(F\langle S\rangle)$, there is a neighborhood $N$ of $T$ such that $B(P', T')$ is $\epsilon$-close to $B(P,T)$ for $T'$ in $N$.  

To show $B(P', T')$ is $\epsilon$-close to $B(P, T)$ in the Gromov-Hausdorff topology, we need to take an arbitrary finite set of points $X$ in $B(P, T)$ and find corresponding points $X'$ in $B(P', T')$ such that the distances between points in $X'$ are within $\epsilon$ of the distances between the corresponding points in $X$.  

For $T'$ close to $T$ in the equivariant Gromov-Hausdorff topology, we can  find points $X'$ in $T'$ with the required properties, but that is not good enough; we need the points $X'$ to be in $B(P',T')$.   We can fix this by projecting each point $x'\in X'$ onto  $B(P',T')$, but we need to be sure that this projection is sufficiently close to $x'$.  

Each projection $p=p(t,s)$ is uniquely determined by the following set of equations (see \cite{GuiLev05}):
\begin{align}\label{eqns}
d(p, sp)&=\ell(s)\notag\\
d(p,tp)&=\ell(t)+2D  \\
d(s^{-1}p,t^{-1}p)&=\ell(s)+\ell(t)+2D \notag \\
d(s^{-1}p, tp)& = \ell(s)+\ell(t)+2D \notag
\end{align}
where $D$ is the distance from $\alpha(s)$ to $\alpha(t)$.  If $T'$ is in the 
$N(\{p\},\{s,t,st^{-1}, st \},\epsilon)$-neighborhood of $T$, then we can find $q'$ in $T'$ satisfying  equations~(\ref{eqns}) up to $\epsilon$. By (\cite{Pau89}), the corresponding lengths $\ell'(s),\ell'(t)$ and $D'$ are within $\epsilon$ of $\ell(s),\ell(t)$ and $D$, so that $q'$ satisfies the analogous equations~(\ref{eqns})$'$ up to $4\epsilon$.  The projection $p'=p'(t,s)$ in $T'$ satisfies the equations~(\ref{eqns})$'$ exactly.  It is an easy exercise to verify that this implies that  $d(p',q')<3\epsilon/{2}$.

Now let $Y=X\cup P$, and let $H=S\cup S^{-1}$.  If $T'$ is in the $N(Y,H,\epsilon)$-neighborhood of $T$, then we can find $q'(t,s)$ and $X'$ in $T'$ with $d(q'(t,s),p'(t,s))<3\epsilon/2$,  $|d(x_1',x_2')-d(x_1,x_2)|<\epsilon$ for all $x_1,x_2\in X$, and $|d(x, p(t,s))-d(x',q'(t,s))| <\epsilon$ for all $x \in X$.  We claim that the projection of  $x'$ onto $B(P',T')$ is within $9\epsilon$ of $x'$ for each $x'\in X'$.

Each point $x$ of $B(P,T)$ is determined by its distances to the points $p(s,t)$ ($x$ is on some straight arc between $p_1=p(t_1,s_1)$ and $p_2=p(t_2,s_2)$; both the fact that it lies on this arc and its position on the arc are determined by its distances to $p_1$ and $p_2$).

If $x$ is on $[p_1,p_2] $ then
\begin{align*}
d(p_1',x') + d( x',p_2') &\leq d(p_1',q_1') + d(q_1',x') + d( x',q_2')+ d( q_2',p_2')\\
&\leq d(p_1,x)+d(x,p_2)+3\epsilon/2+3\epsilon/2 + 2\epsilon\\
&=d(p_1,p_2)+5\epsilon
\end{align*}
and
\begin{align*}
d(p_1',p_2')&\geq d(q_1',q_2')-d(p_1',q_1')-d(p_2',q_2')\\
&\geq d(p_1,p_2)-\epsilon -3\epsilon \\
&=  d(p_1,p_2)-4\epsilon 
\end{align*}
For any three points in a tree, we have $d(a,[b,c])={\frac{1}{2}}(d(b,a)+d(a,c)-d(b,c))$, so $$d(x',[p_1',p_2'])= \frac{d(p_1',x') + d( x',p_2')-d(p_1',p_2')}{2}\leq \frac{9\epsilon}{2}.$$
Since $[p_1',p_2']\subset B(P',T')$, the projection of $x'$ onto $B(P',T')$ is within $9\epsilon/2$ of $x'$, and we may replace $x'$ by this projection.
\end{proof}

Now let $P$ be a set of projections with at least one projection on each axis $\alpha(s)$, for $s\in S$.  We use the convex hull $B(P,T)$ to define a new $F\langle S\rangle$-tree $T_0(P,T)$ as follows. For each $s\in S$, let $b(s)$ be the minimum of the projections $p(t,s)$ in $P$ (where $\alpha(s)$ is oriented in the direction of translation by $s$.)  Form a labeled graph $R(P,T)$ by attaching  an oriented circle to $b(s)$ labeled $s$, whose length  is the translation length of $s$ acting on $T$. We call $R(P,T)$ a  ``stemmed rose."  The subgraph $B(P,T)$ is the  ``stem" and the circles are the  ``petals". Note that for any point $x \in B(P,T)$, there is a canonical identification of $\pi_1(R,x)$ with $F\langle S\rangle$.  Lifting to the universal cover of $R(P,T)$ defines an action $T_0(P,T)$  of  $F\langle S\rangle$ on a tree.   

\begin{lemma}\label{T0} The action $T_0(P,T)$ depends continuously on $T$.
\end{lemma}
\begin{proof}  By Proposition~\ref{BT} the subtree $B(P,T)$ and basepoints $b(s)$ depend continuously on $T$, and by \cite{Pau89} the lengths $\ell(s)$ for $s\in S$ depend continuously on $T$.  Together, this data completely determines $T_0(P,T)$.
\end{proof}  

Now fix a copy of $B(P,T)$ in $T_0(P,T)$.  Then there is a unique equivariant map $f_T: T_0(P,T) \to T$ which is the identity on $B(P,T)$.    This map is a  ``folding" in the sense of Skora: for any segment $[x,y]$ on $T_0(P,T)$, there exist a non-trivial  initial segment $[x,z] \subset [x,y]$, $x \neq z$, such that $f_T$ restricted to $[x,z]$ is an isometry.  We call such a folding map a morphism.
Following \cite{GuiLev05}, we can use this morphism to define a path from $T$ to $T_0(P,T)$ in $\mathcal O(F\langle S\rangle)$ as follows. For $r \in [0, \infty )$, let $T_{r}(P,T)$ be the quotient of $T_0(P,T)$ by the equivalence relation $x \sim_r y$ if $f_T(x)=f_T(y)$ and $f_T[x,y]$ is contained in the ball of radius $r$ around $f_T(x)$.  Since $f_T$ is a morphism, this gives $T_0(P,T)$ for $r=0$.  By \cite{Sko},   $T_{r}(P,T)$ is a tree for all $r$, the action of $F\langle S\rangle$ on $T_0(P,T)$ descends to an action on $T_{r}(P,T)$ and for $r$ sufficiently large, $T_{r}(P,T)=T$.  Furthermore, the path  $p_T: [0,1] \to \mathcal O(F\langle S\rangle)$ defined by $p_T(0)=T$ and $p_T(r)=T_{(1-r)/r}(P,T)$ for $r>0$ is continuous.  

\begin{lemma}\label{unfold} 
The paths $p_T$ define a deformation retraction of $\mathcal O(F\langle S\rangle)$ onto the subspace consisting of universal covers of stemmed roses marked by the generators $S$, which is contractible. 
\end{lemma}

\begin{proof}  
By  Lemma~\ref{T0}, $T_0(P,T)$ depends continuously on $T$.  This implies that the morphisms $f_T\colon T_0(P,T)\to T$ depend continuously on $T$ and then, by Skora's argument (see \cite{GuiLev05}, Proposition 3.4), that the folding paths stay close, i.e. $T_r(P,T)$ is close to $T_r(P',T')$ for all $r$. Therefore these paths give a deformation retraction of all of $\mathcal O(F\langle S\rangle)$ to the subspace of actions covering ``stemmed roses" with petals marked by the generators $S$.  Contracting all stems to points defines a further deformation onto actions covering the standard rose with lengths on its edges.  Such an action is determined by the lengths of the generators, so this space is a product of positive rays, which is contractible.  \end{proof}

We are now ready to produce a contraction  of $\mathcal O(\AG)$.  

By definition of the topology, $\OG$ is homeomorphic to the product $Y \times Q\times \R^\ell $ from Proposition~\ref{topology}.  The space $Q \subset \R^k$ is clearly contractible, so it remains only to show that $Y \times \{0\}$  is contractible where $0$ is the origin in $Q \times \R^\ell  \subset \R^{k+\ell}$.  A point $(T_v)$ in $\prod_{v \in V_0} \mathcal O(F\langle L_v\rangle)$ lies in $Y$ if and only if it satisfies certain compatibility conditions. Namely for each  $w \in V_0$, the minimal $F\langle L\p_w\rangle$-trees in $T_v$ must be equivariantly isometric for all $v \in L_w \cap V_0$.

By Proposition~\ref{topology}, a point  $\X$ in  $\mathcal O(\AG)$ is determined by
\begin{itemize}
\item For each $v\in \G_0$, a tree $T_v$ with $F\langle L_v\rangle$-action and a tree $\TvP$ with $F\langle L\p_v\rangle$-action
\item  For each edge $[v,w]$ in $\G_0$ an $F\langle L_w\p\rangle$-equivariant isometry $i(w,v)$ from $ \TwP$ into $T_v$;
\item For each leaf $w$ of $\G_0$, a real number $\lambda(w)$.
\end{itemize}
If $w$ is cyclic, the isometry $i(w,v)$ is determined by a real number $\Delta_v(w)$. If $w$ is not cyclic, there is a unique equivariant isometry into $T_v$. 
We are assuming $\X$ is in $Y\times\{0\}$, so $\Delta_v(w)=0$ for all $v,w$ and $\lambda(w)=0$ for all leaves $w$.

Given a point $\X$ in $Y\times \{0\}$, we want to produce a new point $\X_0$ in $Y\times \{0\}$
and a morphism $f_{\X}\colon \X_0 \to \X$.  We will show $f_\X$ depends continuously on $\X$, and that the resulting paths $p_\X$ define a deformation retraction of $Y\times \{0\}$ to a contractible subspace.  

{\it Definition of $\X_0$.} Let $v\in V_0$.  We set  $(T_0)_v=T_0(P_v,T_v)$, where $P_v$ is chosen as follows.   If $u\in L_v$ is equivalent to some non-cyclic vertex $w\in \G_0$, include all projections $p(t,u)$ for $t\in [w], t\neq u$.  For all other $u\in L_v$, take all projections $p(t,u)$, for $t\in L_v, t\neq u$. 

If $v$ is not cyclic, set $(T_0)_v\p=T_0(P_v\p,T_v\p)$, where  $P_v\p$ consists of all projections $p(t,u)$ for $t, u \in [v], t\neq u$.  If $v$ is cyclic then $(T_0)_v\p$ is a linear tree, and we fix a basepoint $b(v)$ arbitrarily.  

If $[v,w]$ is an edge of $\G_0$, the image  of $\TvP$ in $T_w$ under $i=i(v,w)$ is the unique $L_v\p$-invariant subtree of $T_v$.  The image of each projection $p(t,u)$ in $P_v\p$ is the analogous projection in $T_v$, so the image of       
 $B(P_v\p,\TvP)$ is contained in $B(P_w,T_w)$, is isometric to the intersection of $B(P_v,T_v)$ with $i(\TvP)$ and has the same basepoints.   We can therefore construct an isometric embedding $i_0(v,w)\colon T_0(P_v\p,T_v\p)\to T_0(P_w,T_w)$.   If $v$ is cyclic, we send the basepoint $b(v)$ of $(T_0)_v\p$ to the basepoint $b_w(v)$ on $\alpha_w(v)$. 
 
The actions $(T_0)_v$ and $(T_0)_v\p$ and isometries $i_0(v,w)$ form a compatible system, giving a point $\X_0$ of $\mathcal O(\AG)$.  The real numbers $\Delta_w(v)$ for this point are all zero, so  in fact $\X_0$ is in $Y\times\{0\}$.  

{\it The morphisms $f_{\X}$ and paths $p_{\X}$.}
For each $v\in V_0$, we have morphisms $f_v\colon T_0(P_v,T_v)\to T_v$ and $f_v\p\colon T_0(P_v\p,T_v\p)\to T_v\p$ as defined above, as well as folding paths $p_v$ from $T_0(P_v,T_v)$ to $T_v$ and  $p_v\p$ from $T_0(P_v\p,T_v\p)$ to $T_v\p$.   All of these objects depend continuously on $\X$ by Lemma~\ref{T0}. 

We next observe that if $v$ and $w$ are connected by an edge, the morphism $f_w\p$ is the restriction of $f_v$ to $T_0(P_w\p,\TwP)$.  Thus points in $T_0(P_w\p,\TwP)$ are identified under $f_w\p$ if and only if they are identified under $f_v$.  Since the image of $\TwP$ is a geodesic subspace of $T_v$,  the folding process produces  equivariant isometries $i_r(w,v)\colon p_w\p(r) \to p_v(r)$ for all $r$, with $\Delta_v(w)=0$ for all cyclic vertices $w$.   

Thus the paths $p_v$ and $p_v\p$ form a compatible system  $p_\X$ of paths in $Y\times \{0\}$ and give a deformation retraction of $Y\times\{0\}$ onto the subspace   whose vertex actions are of the form $T_0(P_v, T_v)$, i.e. with quotient a stemmed rose whose petals are marked by the generators $L_v$.  Simultaneously contracting all stems to a point gives a further deformation retraction onto the subspace $Y_0\times\{0\}$ of $Y\times\{0\}$ whose trees are universal covers of roses (without stems) marked by the generators $L_v$.  These roses are uniquely defined by the translation lengths of the generators, so $Y_0\times\{0\}$ is homeomorphic to a product of positive real rays. It follows that $Y\times\{0\}$ is contractible.  This completes the proof of Theorem~\ref{contractible}.

\bigskip

\section{The action of $Out^0(A_\G)$}


In this section we show that there is a proper (right) action of $Out^0(A_\G)$ on $\OG$.  We continue to assume that $\G$ is connected, triangle-free, and contains more than one edge. 

Let $\phi$ be an element of $Out^0(\AG)$.  Recall from Lemma~\ref{edge-reps} that if
$v$ and $w$ are connected by an edge $e$ in $\G_0$ and $\phi_v$, $\phi_w$ are representatives of $\phi \in Out^0(A_\G)$ preserving $A_{J_v}$ and $A_{J_w}$ respectively, then  there exists $g_v \in A_{J_v}$ and $g_w \in A_{J_w}$ such that $c(g_v)\circ \phi_v = c(g_w) \circ \phi_w$.
 Setting $\phi_e=c(g_v)\circ \phi_v = c(g_w) \circ \phi_w$ gives a representative of $\phi$ which preserves both $A_{J_v}$ and $A_{J_w}$, and hence also $A_{J_e}$.  The automorphism $\phi_e$ is unique up to conjugation by an element of $A_{J_e}$.
 
Now let $\X=\{X_v, X_e, i_{e,v}\}$ be an element of $\OG$ where $i_{e,v}$ denotes the isometric embedding   $X_e \to X_v $ .   Let $X_v^{\phi_v}$ denote the space $X_v$ with the $A_{J_v}$ action twisted by $\phi_v$.  Notice that translation by $g_v$  is an equivariant isometry $t(g_v)\colon X_v^{\phi_v}\to X_v^{\phi_e}$ (where the translation is taken with respect to the original action on $X_v$).  Hence 
\[
 X_e^{\phi_e} \overset {i_{e,v}} \longrightarrow \ X_v^{\phi_e} \overset {t(g_v)^{-1}} \longrightarrow  X_v^{\phi_v}.
 \]
is an equivariant embedding.  We define 
\[
\X \cdot \phi = \{X_v^{\phi_v}, X_e^{\phi_e}, t(g_v)^{-1} i_{e,v} \}.
\]
It is straightforward to check that this is independent of the choices of $\phi_v$ and $\phi_e$:  any other choice gives an equivariantly isometric graph of tree-spaces. 

To see that this defines an action, we must verify that if $\rho$ is another element of $Out^0(A_\G)$, then $(\X \cdot \phi )\cdot \rho  = \X\cdot (\phi\rho) $.   Suppose
\[
\phi_e=c(g_v) \circ \phi_v = c(g_w) \circ \phi_w \qquad  \rho_e=c(h_v) \circ \rho_v = c(h_w) \circ \rho_w.  
\]
Then $(\X \cdot \phi)\cdot \rho = \{(X_v^{\phi_v})^{\rho_v}, (X_e^{\phi_e})^{\rho_e}, t'(h_v)^{-1}t(g_v)^{-1} i_{e,v} \}$ where $t'$ denotes translation with respect to the (twisted) action on $X_v^{\phi_v}$.

To compute $ \X\cdot (\phi\rho) $, note that since $\phi_v$ and $\rho_v$ preserve $A_{J_v}$,  so does their composite, so without loss of generality, we may choose $(\phi\rho)_v=\phi_v \rho_v$.  Observe also that 
\[
\phi_e\rho_e= c(g_v)  \phi_v  c(h_v) \rho_v =c(g_v\phi_v (h_v))  \phi_v \rho_v=c(k_v)(\phi\rho)_v 
 \]
where $k_v= g_v\phi_v (h_v)$ and likewise
\[
\phi_e\rho_e= c(g_w)  \phi_w c(h_w)  \rho_w=c(g_w \phi_v (h_w)) \phi_w \rho_w =c(k_w)(\phi\rho)_w.
 \]
   It follows that we can take
$(\phi\rho)_e=\phi_e \rho_e$ and that $t(k_v)^{-1}=t'(h_v)^{-1}t(g_v)^{-1}$ which completes the argument.

\begin{theorem}\label{proper}  The action of $Out^0(\AG)$ on $\OG$ defined above is proper.
\end{theorem}

\begin{proof}  We must show that for any compact set $C \subset \OG$,
\[
S(C)=\{ \phi \in Out^0(\AG) \mid C \cap C \phi  \neq \emptyset \}
\]
is finite.  Consider first the case in which $\G$ is a single join $\G=U \ast W$.  If $F\langle U \rangle$ and $F\langle W \rangle$ are both non-abelian, then $\OG=\mathcal O(F\langle U \rangle) \times \mathcal O(F\langle W \rangle)$ and the action is proper by \cite{CulVog86}.  If  $U=\{v\}$, then $\OG \cong\R^\ell \times \R_{>0} \times  \mathcal O(F\langle W \rangle)$ 
where the $\R^\ell$ keeps track of the skewing homomorphism $\lambda : F\langle W \rangle \to \R$ and $\R_{>0}$ records the translation length of $v$.  Thus, we can specify an element of $\OG$ by a triple 
$(\lambda,t_v, T_W)$.  The group $Out(\AG)$ decomposes as a semi-direct product, 
 $Out(\AG) \cong \Z^l \rtimes (\Z/2 \times Out(F\langle W \rangle)$.  Denoting an element of $Out(\AG)$ by a triple $(\mathbf z, \epsilon, \phi)$ according to this decomposition, the action of $Out(\AG)$ on $\OG$ is given by
 \[
(\lambda, t_v, T_W) \cdot  (\mathbf z, \epsilon, \phi) = ( \epsilon(\lambda + t_v \mathbf z), t_v, T_W^\phi).
\]
In any compact set in $\OG$, the translation lengths $t_v$ are bounded away from zero. Since the action of  $Out(F\langle W \rangle)$ on $\mathcal O(F\langle W \rangle)$ is proper, it is now easy to see that the action of $Out(\AG)$ on $\OG$ is proper.

Now let $\G$ be arbitrary, and let $R_0$ be the restriction homomorphism from Proposition~\ref{kernel}.  We first show that the kernel $\kerR$ of $R_0$ acts properly on $\OG$.  An element $\phi \in \kerR$ acting on $\X$ fixes all vertex spaces $X_v$ and all edge spaces $X_e$ and acts only on the connecting maps $i_{e,v}$, or equivalently on the factor $Q$ in the product decomposition of $\OG$.  By Remark~\ref{Q}, $Q=\prod_{w \in V_{cyc}} Q_w$, where $Q_w$ can be identified with $\R^{\delta_0(w)}/\R$ with coordinates given by the basepoints $c_v(w)$ for $v \in L_w \cap \G_0$.  We claim that the action of $K_0$ on the $Q_w$ factor is given by the homomorphism $\mu_w: K_0 \to \Z^{\delta_0(w)}/\Z$ described in the proof of Proposition~\ref{kernel}, with $m \in \Z$ acting on $c_v(w)$ as translation by $w^{m}$.  The product $\mu=\prod \mu_w$ is an isomorphism, so it follows easily from the claim that the action of $K_0$ is proper.

To prove the claim, recall from the proof of Proposition~\ref{kernel}, that for $\phi \in K_0$, and $v\in L_w\cap \G_0$, we can write $c(v^n)\circ \phi_v=c(w^m)\circ \phi_w$ for some $n,m$.  The  $v$-factor of 
$\mu_w(\phi)$ is defined to be $-m$.  Rewriting this equation as
\[
c(w^{-m})\circ \phi_v=c(v^{-n})\circ \phi_w
\]
we see that we can choose $\phi_e$ to be $c(w^{-m})\circ \phi_v$.  Then $\phi$ takes the connecting map $i_{e,v} : X_e \to X_v$ to 
\[
t(w^m)\, i_{e,v} = i_{e,v}\, t(w^m) : T_v\p \times T_w\p \to T_v\p \times T_v
\]
and similarly for $i_{e,w}$.  Note that $i_{e,v}\,  t(w^m)$ is the identity on the first factor (assuming that was the case for $i_{e,v}$) and $i_{e,w} t(v^n)$ is the identity on the second. (This was the reason for our particular choice of $\phi_e$.)   Thus, the basepoint $c_v(w)$ in $\X\cdot \phi$ is defined to be the inverse image of the natural basepoint  $b_v(w) \in T_v$ under this connecting map.  It is the  $w^{-m}$ translate of the basepoint $c_v(w)$ in $\X$. This proves the claim.

To show that the whole group $Out^0(\AG)$ acts properly, note  that since the action of $Out(A_{J_v})$ on $\mathcal O(A_{J_v})$ is proper for each $J_v$,   for any compact $C$ the image of $S(C)$ under $R_0$ is finite.  Thus $S(C)$ is contained in a finite set of right cosets $\kerR\phi_1, \dots \kerR\phi_m$.  Let $C_i= C \cup C\phi_i^{-1}$.  If  $\rho \in \kerR$ is such that $C\rho \phi_i  \cap C \neq \emptyset$, then $C_i\rho \cap C_i \neq \emptyset$, so by the paragraph above, there are only finitely many such $\rho$ for each $\phi_i$.  We conclude that $S(C)$ is finite.
\end{proof}


\section{Virtual cohomological dimension}

Since $Out(\AG)$ has torsion, its cohomological dimension is infinte. However, it follows easily from our results that $Out(\AG)$ has torsion-free subgroups of finite index, so that its virtual cohomological dimension (vcd) is defined and finite.  In this section we find upper and lower bounds on this vcd.  

\subsection{The projection homomorphism}

We will use the projection homomorphism 
$$P = \prod_{v \in V_0} P_v : Out^0(\AG) \to \prod Out(F\langle L_v\rangle) $$
defined in Section \ref{sec:kernel}. Recall that $\kerP$ denotes the kernel of $P$ and let $Im(P)$ be the image.

\begin{proposition}\label{tffi} The outer automorphism group of a two-dimensional right-angled Artin group has torsion-free subgroups of finite index.
\end{proposition}

\begin{proof} Since the outer automorphism group of a free group is virtually torsion-free, 
for each $v\in V_0$, we may choose a torsion-free subgroup $H_v$ of finite index in $Out(F\langle L_v\rangle)$.    By Proposition ~\ref{KerP} the kernel $\kerP$ of $P$ is free abelian, so that the preimage $H_0$ of $\prod_v H_v$ in $Out^0(\AG)$ is also torsion-free of finite index.   Since $Out^0(\AG)$ has finite index in $Out(\AG)$, this shows that $Out(\AG)$ itself  is virtually torsion-free. 
\end{proof}

An application of the Hochschild-Serre spectral sequence now gives an upper bound for the virtual cohomological dimension of $Out(\AG)$.

\begin{theorem}\label{vcd-HS}  The virtual cohomological dimension of $Out(\AG)$ satisfies
\[
vcd(Out(\AG)) \leq rank(\kerP) + vcd(Im(P)) \leq  \sum_{v\in V_0} (\cdv -1) ~~+  \sum_{v\in V_0} (2\dv-3).
\]
\end{theorem}
\begin{proof}
Consider the exact sequence
 \[
 1 \to \kerP \to Out^0 (\AG) \to Im(P) \to 1.
 \]
 Restricting this exact sequence to a torsion-free, finite index subgroup of $Out(\AG)$, it follows that for any coefficient module,  the $E^2_{p,q}$-term of the associated Hochschild-Serre spectral sequence is zero for $p > vcd(Im(P))$ or  $q > \text{rank}(\kerP)$. 
The rank of $\kerP$ is   $\sum_{v\in V_0} (\cdv -1)$ by Proposition~\ref{KerP}.   Since $Im(P)$ is a subgroup of  $\prod Out(F\langle L_v\rangle)$ and the vcd of $Out(F_n)$ is equal to $2n-3$ (see \cite{CulVog86}), the vcd of $Im(P)$ is at most $\sum_{v\in V_0} (2\dv-3)$. 
 \end{proof}

If the graph $\G_0$ contains a vertex $v$ which has no leaves attached and is contained in no squares, then the only generators of $Out^0(\AG)$ which affect vertices in $L_v$ are inversions and partial conjuations...there are no transvections onto vertices in $L_v$.  Therefore the image of $Out^0(\AG)$ in $Out(F\langle L_v\rangle)$ is contained in the subgroup $P\Sigma(L_v)$ generated by {\it pure symmetric automorphisms}, i.e. automorphisms which send each generator to a conjugate of itself.  By a result of   Collins \cite{Col89}, the subgroup $P\Sigma(L_v)$ has vcd equal to $\delta(v)-2$.  Thus we can improve the upper bound of Theorem~\ref{vcd-HS} as follows.

\begin{corollary}\label{betterUB}  Let $W_0$ be the vertices of $V_0$ which either have leaves attached or are contained in a square with $v$. The virtual cohomological dimension of $Out(\AG)$ satisfies
\[
vcd(Out(\AG))  \leq   \sum_{v\in V_0}(\cdv + \dv-3) +  \sum_{v\in W_0} (\dv-1).
\]
In particular, if $\G$ has no leaves, triangles or squares, then the virtual cohomological dimension of $Out(\AG)$ satisfies
\[
vcd(Out(\AG)) \leq   \sum_{v\in V_0} (\dv + \cdv  -3).
\]
\end{corollary}

\subsection{Free abelian subgroups}
The rank of a free abelian subgroup of a group gives a lower bound on its virtual cohomological dimension.  We have already exhibited a free abelian subgroup $\kerP$ of $Out(\AG)$, generated by leaf transvections and partial conjugations, but in general this is not the largest one can find.  In this section we exhibit a subgroup which often properly contains $\kerP$.   We begin by identifying three free abelian subgroups of $Aut(\AG)$.

\begin{enumerate}  

\item  {\bf The subgroup A}. Recall that $\G_0$ is a subgraph of $\G$ with one vertex in each maximal equivalence class of vertices, and that the partial order on vertices is given by $v \leq w$ if $lk(v)\subseteq lk(w)$.  For each vertex $v$ which is not in the vertex set $V_0$ of $\G_0$,  choose a vertex $w\in  V_0$ with $v \leq w$, and let $A$ be the free abelian subgroup of $Aut(\AG)$ generated by the left and right transvections $\lambda\colon v\mapsto wv$ and $\rho\colon v\mapsto vw$. The rank of $A$ is   $2|V\backslash V_0|$

\item {\bf The subgroup L}.  Let $ L$ denote the free abelian subgroup of $Aut(\AG)$ generated by leaf transvections.  Then $ L $ has rank $\ell$, the number of leaves of $\G$.

\item {\bf The subgroup C}.  Let $ C$ denote free abelian subgroup of $Aut(\AG)$ generated by partial conjugations by a vertex $v$ of one component of $\G-\{v\}$.  Since this is trivial when a component has only one vertex (which is therefore a leaf),  $ C$ has rank $\sum_{v\in V_0}$ $\cdv-\lv$. 
\end{enumerate}

It is easy to check that all generators of the subgroups $A$, $L$ and $C$ commute and generate a free abelian subgroup of $Aut(\AG)$.  We let $G$ denote the (free abelian) image of this subgroup in $Out(\AG)$, and will now compute the rank of $G$. 
The image of $L$ in $G$ is isomorphic to $L$ and does not intersect the image of the subgroup generated by $A$ and $C$.  The subgroups $A$ and $C$, on the other hand may intersect non-trivially and may contain inner automorphisms.  We introduce the following terminology to keep track of the possibilities:

\begin{notation}  A component of $\G-\{v\}$ is a {\it leaf component} if it contains only one vertex.  It is a {\it twig} if it is not a leaf but is contained in the ball of radius 2 about $v$, and a {\it branch} if it is neither a leaf nor a twig.  Note that if $\G$ is a pentagon, the (unique) component of $\G-\{v\}$ is a branch, since points on the interior of the edge opposite $v$ have distance more than $2$ from $v$.  The number of twigs at $v$ will be denoted $\tv$.  
\end{notation}

\begin{theorem}\label{vcd-lb}If $\G$ is not a star, the subgroup $G$ of $Out(\AG)$ generated by the images of $A$, $L$ and $C$ is free abelian of rank $$2|V\backslash V_0| + \sum_{v\in V_0}(\cdv-\tau(v)-1) $$
  \end{theorem}

\begin{proof} The subgroup of $Aut(\AG)$ generated by $L$ and $A$ has rank $\ell + 2|V\backslash V_0|$, where $\ell=\sum_{v\in V_0} \ell(v)$ is the total number of leaves in $\G$. If $v$ is a separating vertex (which is necessarily in $V_0$), then partial conjugation of a leaf component by $v$ is trivial, and  
partial conjugation of a twig by $v$ is contained in $A$. Partial conjugation of a branch by $v$ is not contained in $A$.    However, the subgroup generated by $A$ and {\it all}  partial conjugations of branches at $v$ contains the inner automorphism associated to $v$.  Thus when we pass to $Out(\AG)$,  partial conjugations at $v$ contribute only $\bv-1$  generators of $G$ which are independent of $A$ and $L$.  
\end{proof}

Theorem \ref{vcd-lb} and Corollary \ref{betterUB}  are summarized in the following corollary.

\begin{corollary}\label{vcd}  The virtual cohomological dimension of $Out(\AG)$ satisfies
\begin{align*}
2|V| + \sum_{v\in V_0}(\cdv-\tau(v)-3) & \leq vcd(Out(\AG)) \\ 
&\leq \sum_{v \in V_0} (\cdv + \dv- 3) + \sum_{v\in W_0}({\dv -1}).
\end{align*}
\end{corollary}

\subsection{Examples}

\begin{example}\label{nmTree} Consider the tree $\G$ in Figure \ref{n,m-tree} consisting of one interior edge with $n$ leaves attached at one vertex $v$ and $m$ leaves attached at the other vertex $w$.  The subtree $\G_0$ is the single interior edge. 
We have $\cdv = \dv = n+1,$ and $\delta_C(w)=\delta(w)=m+1$ and $\tau(v)=\tau(w)=1$ so the left-hand side  of the formula in Corollary~\ref{vcd} is $$2(m+n+2)+(m-3)+(n-3) =3m+3n-2.$$  We have $W_0=V_0$, so the right-hand side is $(2m-1) + (2n-1) + m + n=3m+3n-2.$ Thus, in this example, the upper and lower bounds agree giving a precise computation, $vcd(Out(\AG))=3(n+m)-2$.
  \begin{figure}[ht!] 
\labellist
\small\hair 2pt  
\pinlabel {$n$ leaves} [r] at -15 42
\pinlabel {$m$ leaves} [l] at 232 42
\pinlabel {$v$} at 77 50
\pinlabel {$w$} at 145 50
\endlabellist
\centering
\includegraphics[scale=0.7]{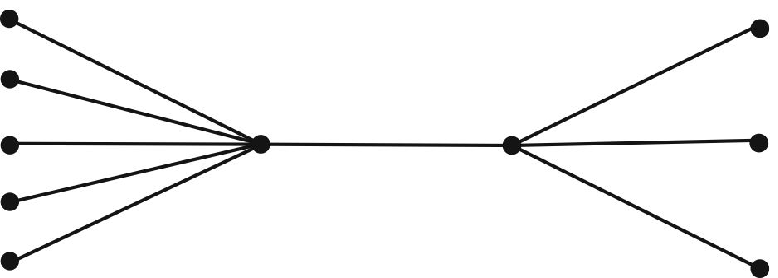}
\caption{} 
\label{n,m-tree}
\end{figure}
\end{example}

\begin{example}
 More generally, suppose that $\G$ is an aribitrary tree. Then $V_0$ is the set of non-leaf vertices of $\G$, $\dv=\cdv$ for all $v$, and there is one twig for each univalent vertex of $\Gamma_0$.  Let $e$ be the number of edges in $\G$, $\ell$ the number of leaves and $\ell_0$ the number of leaves in $\G_0$.   A simple exercise shows that  $e -1=\sum_{v \in V_0} (\dv-1)$; using this, the formulas in Corollary~\ref{vcd} become: 
\[
e-1+2\ell-\ell_0  \leq vcd (Out(\AG)) \leq  e +\ell -3 + \sum_{v\in W_0} (\delta(v)-1).
\]

\end{example}
\begin{example}
 Consider the case of a single join $\G=V \ast W$ with $V=\{v_1, \dots ,v_n\}$ and $W=\{w_1, \dots w_n\}$, $n,m \geq 2$. Then $\G_0$ consists of a single edge from, say, $v_1$ to $w_1$.  The subgroup $K_P$ is trivial, so Theorem~\ref{vcd-HS} gives an upper bound of $(2n-3) + (2m-3)$ on the vcd.  For the lower bound, we note
$\delta_C(v_1)=\delta_C(w_1)= \tau(v_1)=\tau(w_1)=1$, so  the lower bound is $2(m+n)-6$, matching the upper bound.  
\end{example}

\begin{figure}[ht!]
\vskip 15pt
\labellist
\small\hair 2pt  
\pinlabel {$\Gamma$} at 105 -25 
\pinlabel {$\Gamma_0$} at 380 -25
\pinlabel {$V$} at 3 115
\pinlabel {$v$} at -10 58
\pinlabel {$v$} at  290 58 
\pinlabel {$W$} at 147 130
\pinlabel {$W$} at 450 130
\endlabellist
\centering
\includegraphics[scale=0.7]{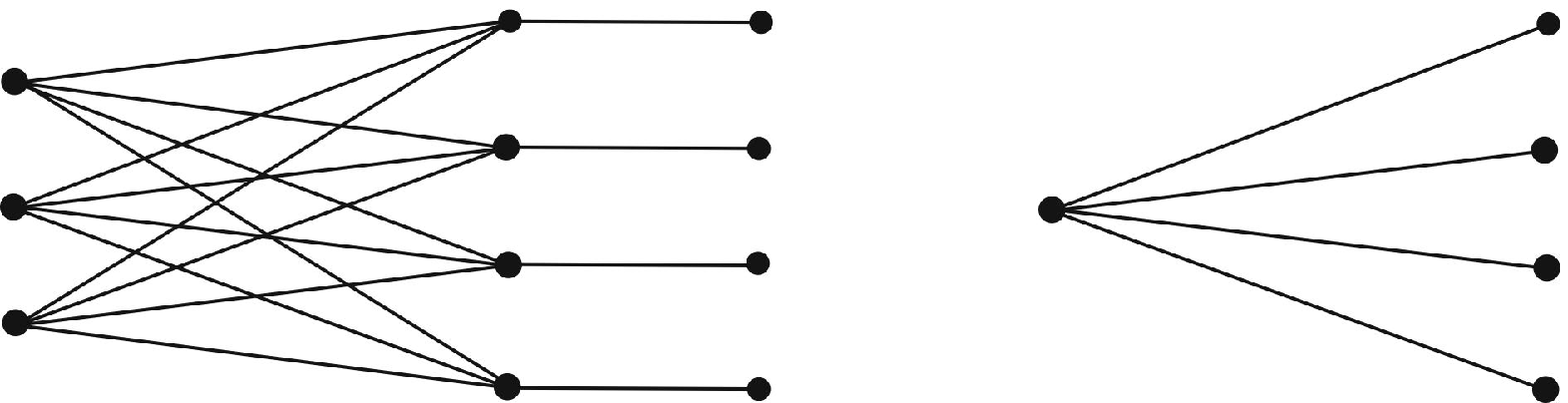}
\vskip 10pt
\caption{} 
\label{small-Z}
\end{figure}

\begin{example}\label{overlap}
 Suppose $\G$ is the graph in Figure \ref{small-Z}, with $n>1$ vertices in $V$ and $m>1$ vertices in $W$.  For $v\in V$ we have $\tau(v)=1$, $\dv=m$ and $\cdv=1$, while for each $w\in W$ we have $\tau(w)=0$, $\delta(w)=n+1$ and $\delta_C(w)=2$.  Thus the rank of $G$ is equal to $2(n+2m) + (-1) + (-m)=2n+3m-3$, giving a lower bound on the vcd, and the upper bound on the vcd is equal to $2mn +2m-3$ since $W_0=V_0$.  So we obtain
\[
3m+2n-3 \leq vcd (Out(\AG)) \leq 2mn +2m-3.
\]
Thus the gap between the upper and lower bounds grows rapidly with $m$ and $n$.
\end{example}

\begin{example}
When the rank of $G$ is equal to the vcd, as in example~\ref{nmTree} above, it follows that $G$ is a maximal rank abelian subgroup in $Out(\AG)$. However, this is not always the case.  For example, suppose $\G$ contains a vertex $v$ such that $\G-\{v\}$ includes a unique leaf $w$ and a large number of non-leaf components.  Then the generating set for $G$ contains 3 tranvections onto $w$ (one leaf-transvection and two non-leaf transvections).  In place of these 3 transvections, one could take all of the partial conjugations by $w$ of a non-leaf component of $\G-\{v\}$.  This makes sense since the non-leaf components of $\G-\{v\}$ are exactly the components of $\G-st(w)$.  One can check that these partial conjugations by $w$ commute with all of the other generators of $G$, giving a larger rank abelian subgroup.  It would be interesting to determine the maximal rank of an abelian subgroup in $Out(\AG)$ and whether that rank is always equal to the virtual cohomological dimension. 
\end{example} 

\subsection{A spine for $\OG$}

The dimension of outer space $\OG$ is in general much larger than the virtual cohomological dimension of $Out(\AG)$.  In the case of a free group $F$, the outer space $\mathcal O(F)$ contains an equivariant deformation retract, called the {\it spine}, with dimension equal to the vcd of $Out(F)$.   In this section we produce a similar spine of $\OG$.  The dimension of this spine is at least as small as the upper bound on the vcd obtained in Theorem~\ref{vcd-HS}, and in several of the examples given in the previous section its dimension is equal to the exact vcd of $Out(\AG)$.

 We begin recalling the construction of  the spine of outer space for a free group.
Since we have not projectivized  $\mathcal O(F)$, it decomposes as a union of open cubes in a cubical complex.  To see this, we view points in $\mathcal O(F)$ as marked, metric graphs, i.e. metric graphs with an isomorphism (determined up to conjugacy)  from $F$ to the fundmental group of the graph.  If $T$ is a metric tree with an $F$-action, then the graph $T/F$ has a natural marking, and the open cube containing this point is paramaterized by  varying the lengths of edges of this graph between 0 and infinity.   Some faces of this cube lie in outer space, others do not.  In particular, if a  face contains a graph with an edge of infinite length, then that face does not lie in $\mathcal O(F)$.  

\begin{remark} Though it plays no role in what follows, we note that the cube complex obtained by including all faces of all cubes is topologically a cone, with cone point the point at which all edges have length zero.  The link of this cone point is the usual simplicial closure of projectivised outer space.
\end{remark}

Let $\bar C$ denote the closure of the open cube $C$ inside of $\mathcal O(F)$.  For two open cubes $C_1,C_2$ in $\mathcal O(F)$, say $C_1 < C_2$ if $C_1$ is a face of $\bar C_2$. The {\it spine} of $\mathcal O(F)$, denoted $Z(F)$, is the simplicial complex whose vertices are labeled by the open cubes in $\mathcal O(F)$ and whose simplices correspond to totally ordered sets of these cubes.   Identifying a vertex $v_C$ in $Z(F)$ with the barycenter of the cube $C$, we can view $Z(F)$ as a subspace of $\mathcal O(F)$. Each open cube $C$ in $\mathcal O(F)$ deformation retracts onto the star of $v_C$ in $Z(F)$ and these retracts fit together to give a retraction of $\mathcal O(F)$ onto $Z(F)$.  Since the action of $Out(F)$ on $\mathcal O(F)$ maps open cubes to open cubes and preserves the partial order, there is an induced action of $Z(F)$.

Now let $\AG$ be an arbitrary right-angled Artin group. Recall from Proposition~\ref{topology} that $\mathcal O(\AG)$  decomposes as a product $\R^\ell\times Q \times Y$  where $Y$ is a subspace of the product  $\prod_{v \in V_0} \mathcal O(F\langle L_v\rangle)$ of outer spaces for the free groups $F\langle L_v\rangle$.  Since a product of cubes is cube, this product of outer spaces is as a union of open cubes, where a cube $C=\prod C_v$ corresponds to a specified marked graph for each $F\langle L_v\rangle$, and the edge lengths give coordinates for the cubes.
 
\begin{lemma}  The intersection of each closed cube $\bar C=\prod \bar C_v$ with $Y$ is a convex cell.
\end{lemma}  
\begin{proof}   A collection of trees $\{T_v\}$ lies in $Y$ if and only if the minimal $F\langle L_w\p\rangle$-subtrees are equivariantly isometric in all $T_v$ with $v \in L_w$.   By  \cite{KrsVog93}, each of these subtrees is uniquely determined by the translation lengths of a finite set of elements of $F\langle L_w\p\rangle$.  For each $T_v$, these translation lengths are given by a linear combination of the edge lengths of the graph $T_v/F\langle L_v\rangle$.  Since the edge lengths give coordinates for the cube $\bar C$,  the intersection of $Y$ with $\bar C$ is  given by a finite set of linear equalities.
\end{proof}

It follows from the proof above that if $C_1$ is a face of $\bar C_2$ and their intersection with $Y$ is non-empty, then $Y \cap C_1$ is a face of $Y \cap \bar C_2$.

\begin{definition} We define the {\it spine of $Y$} to be the geometric realization of the poset of cells $Y\cap C$ partially ordered by the face relation.  We denote this spine by $Z(\AG)$.   
\end{definition}

\begin{proposition} \label{spine}
The action of $Out^0(\AG)$ on $\OG$ descends to a proper action of $Im(P)$ on $Y.$  With respect to this action, 
 $Z(\AG)$ can be identified with a piecewise linear  $Im(P)$-invariant subspace of $Y.$ This subspace is a deformation retract of $Y$ hence,  in particular,  $Z(\AG)$ is contractible.
\end{proposition}

\begin{proof}  Let $\X=\{X_v,X_e,i_{e,v}\}$  be a point in $\OG$ and 
$\X \cdot \phi = \{X_v^{\phi_v}, X_e^{\phi_e}, t(g_v)^{-1} i_{e,v} \}$ its translate by $\phi$.  Recall that 
$X_v$ is a product of trees $T_v \times T_v^\perp$ with an action of $F\langle L_v\rangle \times F(L_v^\perp)$.   Though this action is not necessarily a product action, it projects to an action of $F\langle L_v\rangle$ on $T_v$.
The projection of $\X$ on $Y$ is given by the resulting set of actions $\{T_v\}$.  The twisted tree-space $X_v^{\phi_v}$ is a product of the same two underlying trees with the action twisted by  $\phi_v$. The new action of $F\langle L_v\rangle$ on $T_v$ depends only on the projection of $\phi_v$ to $Out(F\langle L_v\rangle)$.
The first statement of the lemma follows.

 The action of $Out(F\langle L_v\rangle)$ on $\mathcal O(F\langle L_v\rangle)$ is cellular and the stabilizer of any cell is finite.  Hence the same is true of the action of  $Im(P)$ on $\prod_{v \in V_0} \mathcal O(F\langle L_v\rangle)$.  By the discussion above, $Im(P)$ preserves $Y$ and  hence it takes open cells $C \cap Y$ to open cells and preserves the face relation. 

Let $p$ be a point of $C\cap Y$.  The orbit of $p$ intersects $C\cap Y$ in a finite set of points.  Since $C\cap Y$ is convex, the barycenter $\bar p$ is a point of $C\cap Y$ which is invariant under the stabilizer of $C\cap Y$, and the entire orbit of $\bar p$ intersects each cell $C'\cap Y$ in at most point.   
  It follows that we can chose one point $x_C$ of each cell $C \cap Y$ such that the set of points $\{x_C\}$ is $Im(P)$-invariant.  Now identify these points with the vertices of $Z(\AG)$ in the obvious way.  Then for any simplex $\sigma$ of $Z(\AG)$, the  vertices of $\sigma$ lie in the closure of a single cell $C \cap Y$ and their linear span forms a simplex in $Y$. 
 The resulting simplicial complex is isomorphic to $Z(\AG)$.  As in the case of the spine for a free group, retracting each cell $C \cap Y$ linearly onto the star of the vertex $x_C$ gives a deformation retraction of $Y$ onto $Z(\AG)$.
 \end{proof}
 
\begin{proposition}\label{vcd-ub}  The virtual cohomological dimension of $Out(\AG)$ satisfies
\[
vcd(Out(\AG)) \leq rank(\kerP) + dim ~Z(\AG) 
\]
\end{proposition}

 \begin{proof} It follows immediately from Proposition~\ref{spine} that $Im(P)$ has vcd bounded by the dimension of $Z(\AG)$.  The result now follows from Theorem~\ref{vcd-HS}.
    \end{proof}

 \begin{proposition}  The dimension of $Z(\AG)$ is at most $\sum_{v}2\dv-3$.
 \end{proposition}
 
 \begin{proof}  The dimension of a cube $C=\prod C_v$ is at least $\sum_v \dv$ and at most  $\sum_v 3\dv-3$ since each cube $C_v$ has dimension at least $\dv$ and at most $3\dv-3$.  Therefore the longest possible chain of inclusions of cells $Y\cap C$ is $(\sum_v 2\dv-3) + 1$, so that the dimension of $Z(\G)$ is at most $\sum_v 2\dv-3$.
  \end{proof}
  
In fact,  $Z(\AG)$ is naturally isomorphic to a subcomplex of a simplicial subdivision of the product $\prod_v Z(F\langle L_v\rangle)$ of spines for the outer spaces associated to the vertices $v$ of $\G_0$.  If the links of vertices of $\G_0$ have large overlap, as in Example~\ref{overlap}, $Z(\AG)$ will be much smaller than the full product, though it can be shown in this example that they have the same dimension.

On the other hand, in  the case that $\G$ is a tree, we claim that $Z(\AG) =\prod Z(F\langle L_v\rangle)$.  To verify this claim, we must show that $Y$ intersects every open cell in $\prod \mathcal O(F\langle L_v\rangle)$.  
For $\G$ a tree, two links $L_v$ and $L_w$ intersect either in exactly one point (if $v,w$ are distance 2 apart) or not at all. If $L_v \cap L_w=\{u\}$, we will say that $T_v$ and $T_w$ are \emph{compatible} if the translation lengths of $u$ in $T_v$ and $T_w$ agree.  A point in $Y$ is a $V_0$-tuple $(T_v)$ of compatible trees. 

A cell in $\mathcal O(F\langle L_v\rangle)$ is invariant under scaling, i.e., if $T_v$ lies in an open cell $C_v$, then so does the tree obtained by scaling the metric on $T_v$ by any $r>0$.  Thus it suffices to show that  any $V_0$-tuple of trees $(T_v)$ can be made compatible by rescaling.   To do this, fix a pair of adjacent vertices $v,u$ in $\G_0$. Every vertex in $\G_0$ is even distance from exactly one of these two vertices.  If $w$ is distance $2n$ from $v$, then there is a unique sequence of vertices $v=w_0, w_1, \dots w_n=w$ such that the link of $w_{i-1}$ intersects the link of $w_i$ in a vertex.  Starting with $T_v$, we can inductively scale each  $T_{w_i}$ to be compatible with the previous one.  Similarly, for vertices  at distance $2n$ from $u$.  The resulting collection of trees defines a point  in $Y$.  

Note that once the metrics on $T_u$ and $T_v$ are fixed, the scaling on the remaining trees is uniquely determined.  Thus, modulo scaling the two base trees, a point in $Y$ corresponds to a point in the product of the \emph{projectivized} outer spaces $\overline{ \mathcal O}(F\langle L_v\rangle)$.
Summarizing, we have shown

\begin{corollary}\label{tree-spine}  If $\G$ is a tree, then $Y \cong \R^2\times \prod \overline{\mathcal O}(F\langle L_v\rangle)$ and $Z(\AG) \cong \prod Z (F\langle L_v\rangle)$, where the products are taken over the non-leaf vertices in $\G$. 
 \end{corollary}

\medskip\noindent{\bf Acknowledgments:}  Ruth Charney was partially supported by NSF grant
DMS-0405623. John Crisp was partially supported by grant ACI JC1041 from the Minist\`ere D\'el\'egu\'e \`a la Recherche et aux Nouvelles Technologies de la France. 
Karen Vogtmann was partially supported by NSF grant DMS-0204185.

\def\cprime{$\prime$}

\end{document}